\DeclareSymbolFont{AMSb}{U}{msb}{m}{n}
\documentclass[reqno,centertags,11pt,a4paper,noamsfonts]{amsart}
\pdfoutput=1
\usepackage[svgnames]{xcolor}
\usepackage{graphicx}
\usepackage[english]{babel}
\usepackage[babel=true,expansion=alltext,protrusion=alltext-nott,final]{microtype}
\usepackage[margin={3cm},marginparwidth={2.5cm}]{geometry}
\usepackage[utf8]{inputenc}
\usepackage{textcomp}
\usepackage[sb]{libertine}
\usepackage[libertine,bigdelims,vvarbb]{newtxmath}
\useosf
\usepackage[varqu,varl]{inconsolata}
\usepackage[supstfm=libertinesups,%
       supscaled=1.2,%
       raised=-.13em]{superiors}
\usepackage{mathtools}
\usepackage[T1]{fontenc}
\usepackage{textcomp}
\usepackage{amsthm}
\usepackage[inline]{enumitem}
\numberwithin{equation}{section}
\usepackage[normalem]{ulem}

\usepackage{float}  
\usepackage{caption}
\usepackage{subcaption} 
\usepackage{xspace}         
\usepackage[scientific-notation=true]{siunitx}      
\usepackage{pstricks}
\usepackage{tikz}
\usetikzlibrary{positioning} 
\usetikzlibrary{calc}
\usepackage{pgfplots}
\pgfplotsset{width=10cm,compat=1.9}
\usepackage{bm}
\usepackage{sidecap}
\usepackage{graphicx,color}

\usepackage{amsmath}
\DeclareFontFamily{U}{mathx}{}
\DeclareFontShape{U}{mathx}{m}{n}{<-> mathx10}{}
\DeclareSymbolFont{mathx}{U}{mathx}{m}{n}
\DeclareMathAccent{\widehat}{0}{mathx}{"70}
\DeclareMathAccent{\widecheck}{0}{mathx}{"71}

\providecommand{\mr}[1]{\href{http://www.ams.org/mathscinet-getitem?mr=#1}{MR~#1}}
\providecommand{\zbl}[1]{\href{https://zbmath.org/?q=an:#1}{Zbl~#1}}

\providecommand{\doi}[2]{\href{https://doi.org/#1}{#2}}

\newcommand{\CC}{\mathbb{C}}

\newcommand{\C}{\mathcal{C}}

\newcommand{\ii}{\imath}

\newcommand{\e}{\mathbf{e}}

\definecolor{light_gray}{gray}{0.75}
\definecolor{lighter_gray}{gray}{0.5}
\colorlet{light_blue}{blue!20}
\definecolor{dark_green}{rgb}{0.0, 0.6, 0.0}
\definecolor{royal_blue}{rgb}{0.0, 0.22, 0.66}
\definecolor{salmon}{rgb}{1.0, 0.55, 0.41}
\definecolor{gold}{rgb}{0.8, 0.63, 0.21}
\definecolor{navy_blue}{rgb}{0.0, 0.0, 0.5}
\definecolor{crimson}{rgb}{0.79, 0.0, 0.09}
\definecolor{amethyst}{rgb}{0.6, 0.4, 0.8}
\definecolor{alizarin}{rgb}{0.82, 0.1, 0.26}
\definecolor{amaranth}{rgb}{0.9, 0.17, 0.31}
\definecolor{azure}{rgb}{0.0, 0.5, 1.0}
\definecolor{canaryyellow}{rgb}{0.82, 0.41, 0.12}
\definecolor{carrotorange}{rgb}{0.8, 0.33, 0.0}
\definecolor{cadmiumgreen}{rgb}{0.0, 0.42, 0.24}
\definecolor{copper}{rgb}{0.72, 0.45, 0.2}
\definecolor{aqua}{rgb}{0.5, 1.0, 0.83}
\definecolor{awesome}{rgb}{1.0, 0.13, 0.32}
\definecolor{candyapplered}{rgb}{1.0, 0.03, 0.0}
\definecolor{caribbeangreen}{rgb}{0.0, 0.8, 0.6}
\definecolor{indigo}{rgb}{0.0, 0.25, 0.42}

\DeclareMathOperator{\weaklystar}{\rightharpoonup\kern-2.2ex ^* \, \,}

\def\XXint#1#2#3{{\setbox0=\hbox{$#1{#2#3}{\int}$ }
\vcenter{\hbox{$#2#3$ }}\kern-.6\wd0}}

\newcommand{\R}{\mathbb R}
\newcommand{\N}{\mathbb N}
\newcommand{\Z}{\mathbb Z}
\renewcommand{\C}{\mathbb C}

\newcommand\norm[1]{\lVert #1 \rVert}

\newcommand\inner[1]{\langle #1 \rangle}

\newcommand{\ra}{\rightarrow}

\newcommand{\mL}{\mathrm{L}}
\newcommand{\re}{\mathrm{Re}\,}
\renewcommand{\phi}{\varphi}
\newcommand{\dom}{\mathrm{dom}\,}

\newcommand{\ee}{\mathrm{e}}

\theoremstyle{plain}
\newtheorem{theorem}{Theorem}[section]
\newtheorem{proposition}[theorem]{Proposition}
\newtheorem{corollary}[theorem]{Corollary}
\newtheorem{lemma}[theorem]{Lemma}

\newtheorem*{theorem*}{Theorem}

\theoremstyle{definition}
\newtheorem{definition}[theorem]{Definition}
\newtheorem{remark}[theorem]{Remark}
\newtheorem*{remark*}{Remark}

\usepackage[pagebackref=false]{hyperref}
\hypersetup{
	linktoc=all,
	bookmarksdepth=3,
	colorlinks=true,
	linkcolor=Green,citecolor=Orange,urlcolor=DarkBlue,
	pdftitle={Generalized Hadamard three lines lemma},
	pdfauthor={Thiago Carvalho Corso}, 
	pdfsubject={},
	pdfkeywords={}} 
\begin{document}
\numberwithin{table}{section}
\title{A generalized three lines lemma in Hardy-like spaces}

\author[T.~Carvalho~Corso]{Thiago Carvalho Corso}
\address[T.~Carvalho Corso]{Institute of Applied Analysis and Numerical Simulation, University of Stuttgart, Pfaffenwaldring 57, 70569 Stuttgart, Germany}
\email{thiago.carvalho-corso@mathematik.uni-stuttgart.de}

\date{\today}
\thanks{\emph{Funding information}: DFG -- Project-ID 442047500 -- SFB 1481.\\[1ex]
\textcopyright 2024 by the authors. Faithful reproduction of this article, in its entirety, by any means is permitted for noncommercial purposes.}
\begin{abstract}
In this paper we address the following question: given a holomorphic function with prescribed $L^p(\R)$ and $L^q(\R)$ norm (with $1\leq p,q \leq \infty$) along two parallel lines in the complex plane, then what is the maximum value that this function can achieve at a given point between these lines. Here we show that this problem is well-posed in suitable Hardy-like spaces on the strip. Moreover, in this setting we completely solve this problem by providing not only an explicit formula for the optimizers but also for the optimal values. In addition, we briefly discuss some applications of these results to interpolation theory and to Lieb-Thirring inequalities.
\end{abstract}
\maketitle
\setcounter{secnumdepth}{3}
\section{Introduction}

In this paper we study the following question: let $S \subset \C$ be an open strip between two parallel lines $L_1 , L_2 \subset \C$ in the complex plane and let $1\leq p,q\leq \infty$ and $a,b >0$, then what is the maximum value that a holomorphic function $h: \overline{S} \rightarrow \C$ satisfying the constraints
\begin{align*}
    \norm{h}_{L^p(L_1)}^p = \int_{L_1} |h(z)|^p  \mathcal{H}^1(\mathrm{d}z) = a^p\quad\mbox{and} \quad \norm{h}_{L^q(L_2)}^q = \int_{L_2} |h(z)|^q  \mathcal{H}^1(\mathrm{d}z) = b^q,
\end{align*}
where $\mathcal{H}^1$ denotes the one-dimensional Hausdorff measure, 
can achieve at the point $z_0 \in S$. In other words, what is the value
\begin{align}
    M_{a,b}^{p,q}(z_0;S) \coloneqq \sup \{ |h(z_0)| : \norm{h}_{L^p(L_1)} = a \quad\mbox{and}\quad \norm{h}_{L^q(L_2)}= b\} \label{eq:supproblem}
\end{align}
where the supremum is taken over a suitable set of holomorphic functions $h: \overline{S} \rightarrow \C$.

In fact, after translating, rotating, and re-scaling our original strip (and re-scaling $a,b>0$ accordingly), we can assume that $L_1= \R$, $L_2 = \R +\ii$, and $S$ is the horizontal strip
\begin{align*}
    S\coloneqq \{z = x + \ii y \in \C: 0<y <1\}. 
\end{align*}
Moreover, if we consider the supremum over any space of holomorphic functions which is invariant under the transformation 
\begin{align}
    h \mapsto h(z-\tau) \ee^{\ii \lambda z + \beta}\quad \mbox{ for any $\lambda \in \R, \beta \in \C$ and $\tau \in \R$,} \label{eq:symmetries}
\end{align}
then problem~\eqref{eq:supproblem} (for all $a,b>0$) with $z_0 = x + \ii \alpha$ for some $0<\alpha<1$ is equivalent\footnote{Indeed, for any $h$ in a such a space, we can find $\lambda,\beta, \tau \in \R$ such that $h^{\lambda,\beta,\tau}(z) = h(z-\tau) \ee^{\ii \lambda z + \beta}$ satisfy the required $L^p$ and $L^q$ constraints along the lines $\R$ and $\R+\ii$ and such that $|h^{\lambda,\beta,\tau}(x+\ii \alpha)| \approx \sup_{x\in \R} |h^{\lambda,\beta,\tau}(x+\ii \alpha)| = \norm{h^{\lambda,\beta,\tau}_\alpha}_{L^\infty(\R)}$, which implies
\begin{align*}
    M_{a,b}^{p,q}(x+\ii \alpha;S) = a^{1-\alpha} b^\alpha H_{p,q}(\alpha), \quad \mbox{where $M_{a,b}^{p,q}$ and $H_{p,q}$ given respectively by \eqref{eq:supproblem} and~\eqref{eq:pq3lineproblem0}.}
\end{align*}} to the problem
\begin{align}
     \sup_{h} \frac{\norm{h_\alpha}_{L^\infty(\R)}}{\norm{h_0}_{L^p(\R)}^{1-\alpha} \norm{h_1}_{L^q(\R)}^\alpha},\quad  1\leq p,q \leq \infty, \label{eq:pq3lineproblem0}
\end{align}
where $h_y$ denotes the function $x\mapsto h_y(x) = h(x+\ii y)$. It turns out, however, that the supremum in~\eqref{eq:pq3lineproblem0}, when taken over the space of holomorphic functions in a neighborhood of $\overline{S}$, is not bounded. This can be seen, for instance, by considering the function
\begin{align*}
    h(z) = \exp(-\ii \exp(\pi z)) \frac{1}{z^2+4}, \quad z\in S.
\end{align*}

Fortunately, in the case $p=q=\infty$, there is a natural class of functions where the above problem is well-posed and has a simple solution. Precisely, if we consider the Hardy space 
\begin{align*}
    \mathbb{H}^\infty(S) \coloneqq \{h: S \rightarrow \C \mbox{ holomorphic with } \sup_{z \in S} |h(z)| < \infty\},
\end{align*}
then 
\begin{align}
    H_{\infty,\infty}(\alpha) \coloneqq \max_{h \in \mathbb{H}^\infty(S)} \frac{ \norm{h_{\alpha}}_{L^\infty(\R)}}{\norm{h_0}_{L^\infty(\R)}^{1-\alpha} \norm{h_1}_{L^\infty(\R)}^\alpha} = 1, \quad \mbox{for any $0<\alpha<1$,}  \label{eq:hadamard}
\end{align}
where the boundary values $h_0$ and $h_1$ are well-defined in a weak sense (to be clarified later). Moreover, the unique maximizer up to the transformation in~\eqref{eq:symmetries} is the constant function $h(z) =1$. This result is known as Hadamard's three-lines lemma (or theorem) \cite[Theorem 12.8]{Rud87}, and is a classical result in complex analysis that plays a fundamental role in (complex) interpolation theory (see, e.g., \cite[Appendix to IX.4]{RS75} and \cite[Section 1.3.2]{Gra14}).  

Recently, it has been shown by the author and Tobias Ried \cite{CCR24} that the variational problem
\begin{align}
    H_{\infty,2}(\alpha) \coloneqq \sup_{h\in \mathbb{H}^{\infty,2}(S)} \frac{ \norm{h_{\alpha}}_{L^\infty(\R)}}{\norm{h_0}_{L^\infty(\R)}^{1-\alpha} \norm{h_1}_{L^2(\R)}^\alpha}, \label{eq:problem2}
\end{align}
where $\mathbb{H}^{\infty,2}(S)$ denotes a certain Hardy-like subset of holomorphic functions on the strip (see definition below), also admits an unique (up to symmetries) maximizer with an explicit analytic formula. In this case, however, the maximizer has a non-trivial dependence on $\alpha \in (0,1)$ and the formula is not so simple to state. Moreover, despite the similarities between problems~\eqref{eq:hadamard} and \eqref{eq:problem2}, the proof of the results for the latter is completely orthogonal to the standard proof of Hadamard's three lines lemma, which relies on the classical maximum principle for bounded holomorphic functions (or the Phragm\'en-Lindel\"of principle \cite[Theorem 6.1]{Lan99}).

In this paper, our goal is to extend the approach developed in \cite{CCR24} to deal with the general class of variational problems
\begin{align}
	H_{p,q}(\alpha) \coloneqq \sup_{h \in \mathbb{H}^{p,q}(S)} \frac{\norm{h_\alpha}_{L^\infty}}{\norm{h_0}_{L^p(\R)}^{1-\alpha} \norm{h_1}_{L^q(\R)}^\alpha}, \quad 0<\alpha <1,\quad  1\leq p,q \leq \infty, \label{eq:pq3lineproblem}
\end{align}
where $\mathbb{H}^{p,q}(S)$ is the following natural generalization of the (translation invariant) Hardy spaces on the strip.
\begin{definition}[$\mathbb{H}^{p,q}(S)$ spaces] \label{def:Hardy0}
Let $1\leq p,q \leq \infty$, we denote by $\mathbb{H}^{p,q}(S)$ the space of holomorphic functions $h:S \rightarrow \C$ satisfying the bound
\begin{align*}
\norm{h}_{\mathbb{H}^{p,q}} \coloneqq \sup_{0<y<1} \inf_{f+g = h_y} \{(1-y)^{-1} \norm{f}_{L^p(\R)} + y^{-1} \norm{g}_{L^q(\R)}\} <\infty . 
\end{align*}
\end{definition}
As we shall show later, these are the correct spaces to study problem~\eqref{eq:pq3lineproblem}. Moreover, the main contributions of this paper can be summarized as follows:
\begin{enumerate}
    \item We introduce and prove some of the fundamental properties of the Hardy-like spaces $\mathbb{H}^{p,q}(S)$. In particular, we provide a characterization of its boundary values and prove a Poisson representation formula analogous to the one for the classical Hardy spaces on the half-plane (see Section~\ref{sec:Hardy}).
    \item We completely solve problem~\eqref{eq:pq3lineproblem} by providing not only an analytic formula for the optimizers (Theorem~\ref{thm:optimizers}) but also for the optimal values (Theorem~\ref{thm:optimum} and Remark~\ref{rem:analyticformula}). 
    \item We briefly discuss some applications of these results to interpolation theory and to spectral inequalities such as the Lieb-Thirring and Cwikel-Lieb-Rozenblum inequalities.
\end{enumerate}

\subsection{Main results} Let us now state the main results of this paper precisely. We start with two relations between the optimal values $H_{p,q}(\alpha)$ for different values of $p,q\in [1,\infty]$. The first relation is a duality relation between $H_{p,q}(\alpha)$ and $H_{p^\ast,q^\ast}(\alpha)$, where $p^\ast, q^\ast \in [1,\infty]$ denote the H\"older conjugate exponents of $p$ and $q$.

\begin{theorem}[Duality relation] \label{thm:duality} Let $\alpha \in(0,1)$ and $1\leq p,q \leq \infty$, then the variational problem defined in \eqref{eq:pq3lineproblem} admits an unique optimizer $h \in \mathbb{H}^{p,q}(S)$ up to the transformation $h(z) \mapsto h(z-\tau) \ee^{\ii \lambda z + \beta}$ with $\tau, \lambda \in \R$ and $\beta \in \C$. Moreover, we have the relation 
\begin{align}
    \label{eq:dualityrelation} H_{p,q}(\alpha) H_{p^\ast,q^\ast}(\alpha) = \frac{1}{4 \sin(\pi \alpha) \alpha^\alpha (1-\alpha)^{1-\alpha}},
\end{align}
where $1\leq p^\ast, q^\ast \leq \infty$ denote the H\"older conjugate exponents of $p$ and $q$, i.e., the unique exponents satisfying $\frac{1}{p^\ast} + \frac{1}{p} = 1 = \frac{1}{q^\ast} + \frac{1}{q}$. 
\end{theorem}

An interesting consequence of Theorem~\ref{thm:duality} is that it allows us to explicitly compute the values of $H_{p,q}(\alpha)$ in some specific cases. Indeed, from Hadamard's three lines lemma and~\eqref{eq:dualityrelation} we have
\begin{align*}
    H_{\infty,\infty}(\alpha)=1, \quad H_{1,1}(\alpha) = \frac{1}{4 \sin(\pi \alpha) \alpha^\alpha (1-\alpha)^{1-\alpha}}, \quad \mbox{and}\quad  H_{2,2}(\alpha) = \sqrt{\frac{1}{4 \sin(\pi \alpha) \alpha^\alpha (1-\alpha)^{1-\alpha}}}.
\end{align*}
Similarly, from the relation $H_{p,q}(\alpha) = H_{q,p}(1-\alpha)$ (which follows by flipping the strip), we obtain
\begin{align*}
   H_{p,p^\ast}(1/2) = = \frac{1}{\sqrt{2}}, \quad \mbox{for any $1\leq p \leq \infty$.}
\end{align*}

For the general case $1\leq p,q\leq \infty$, however, the value of $H_{p,q}(\alpha)$ can not be obtained from Theorem~\ref{thm:duality}. Nevertheless, in this case we have the following result.
\begin{theorem}[Optimal value] \label{thm:optimum} For any $\alpha \in (0,1)$, let
\begin{align}
    \mathscr{P}_\alpha(x) = \frac{1}{2} \frac{\sin(\pi \alpha)}{\cosh(\pi x) - \cos(\pi \alpha)} \label{eq:Poissondef}
\end{align}
be the Poisson kernel on the strip $S$. Then we have
\begin{align}
    \log H_{1,\infty}(\alpha) =  \int_\R  \mathscr{P}_{\alpha}(x) \log \left(\frac{\mathscr{P}_{\alpha}(x)}{1-\alpha}\right) \mathrm{d} x.  \label{eq:H1infvalue}
\end{align}
Moreover, for any $1\leq p,q\leq \infty$ we have
\begin{align*}
    \log H_{p,q}(\alpha) = \frac{1}{p} \int_\R \mathscr{P}_{\alpha}(x) \log \left(\frac{\mathscr{P}_{\alpha}(x)}{1-\alpha}\right) \mathrm{d} x + \frac{1}{q} \int_\R \mathscr{P}_{1-\alpha}(x) \log \left(\frac{\mathscr{P}_{1-\alpha}(x)}{\alpha}\right) \mathrm{d} x .
\end{align*}
\end{theorem}

\begin{remark}[Analytic formula with dilogarithm function] \label{rem:analyticformula} The integral in \eqref{eq:H1infvalue} can not be evaluated in terms of elementary functions. Nevertheless, it can be expressed in terms of the dilogarithm function
\begin{align*}
    \mathrm{Li}_2(z) \coloneqq - \int_0^z \frac{\log(1-u)}{u} \mathrm{d} u, \quad z\in \C\setminus [1,\infty)
\end{align*}
as follows (for a proof, see Appendix~\ref{app:analyticformula}) : 
\begin{align}
    \log H_{1,\infty}(\alpha) = -(1-\alpha) \log\left(4(1-\alpha) \sin(\pi \alpha)\right) + \frac{1}{\pi} \mathrm{Im}\, \mathrm{Li}_2(\ee^{\ii 2\pi \alpha}) \label{eq:analyticformula} 
\end{align}
In particular, since only the imaginary part of $\mathrm{Li}_2$ appears in the formula above, we have
\begin{align}
    H_{\infty,2}(\alpha) = \frac{1}{(4 \alpha \sin(\pi \alpha))^{\frac{\alpha}{2}}} \exp\left(\frac{1}{2\pi} \mathrm{CI}_2\left(2\pi(1-\alpha)\right)\right), \label{eq:analyticformula2case}
\end{align}
where $\mathrm{CI}_2$ denotes the Clausen integral (or function) \cite[Chapter 4]{Lew81},
\begin{align}
    \mathrm{CI}_2(2\pi \alpha) \coloneqq - \int_0^{2\pi \alpha} \log(2|\sin(u/2)|) \mathrm{d}u. \label{eq:clausendef}
\end{align}
\end{remark}

Our last main result is an analytic formula for the optimizers of problem~\eqref{eq:pq3lineproblem}.
\begin{theorem}[Optimizers] \label{thm:optimizers} Let $\alpha \in (0,1)$ and $1\leq p, q\leq \infty$. Then the unique optimizer of \eqref{eq:pq3lineproblem} up to the transformation $h(z) \mapsto h(z-\omega) \ee^{\ii \lambda z + \beta}$ with $\omega, \lambda \in \R$, $\beta \in \C$, is given by
\begin{align*}
    h(z) = \ee^{\phi_{\alpha,p,q}(z)}, \quad z\in S,
\end{align*}
where $\phi_{\alpha,p,q}: S \rightarrow \C$ is the unique holomorphic function satisfying $\mathrm{Im}\, \phi_{\alpha,p,q}(\ii \alpha) =0$ and
\begin{align*}
    \re \phi_{\alpha,p,q}(x+\ii y) = \frac{1}{p} \mathscr{P}_y \ast \log\left(\frac{\mathscr{P}_\alpha}{1-\alpha}\right)(x) + \frac{1}{q} \mathscr{P}_{1-y} \ast \log\left(\frac{\mathscr{P}_{1-\alpha}}{\alpha}\right) (x), \quad x+\ii y \in S, 
\end{align*}
where $\mathscr{P}_y$ is the Poisson kernel defined in \eqref{eq:Poissondef}, and $\ast$ denotes the standard convolution in $\R$.
\end{theorem}

\begin{remark*} In the case $p=\infty$ and $q=2$, an alternative formula for the optimizers, which involves a Blaschke factor on the upper-half plane and the Fourier transform of a principal-value like distribution, can be found in \cite[Theorem 1.3]{CCR24}. 
\end{remark*}
\subsection{Applications}

In this section, we briefly mention some applications of our main results. 

The first application we present is a sharp weighted version of inequality~\eqref{eq:pq3lineproblem}. Note that, unlike the $\mL^\infty$ norm, the $\mL^p$ norms (for $1\leq p <\infty$) along the boundary are not conformally invariant. Therefore, the following weighted inequality could be useful as it allows for some degree of conformal invariance in the general three-lines lemma. 
\begin{theorem}[Weighted three-lines inequality] \label{thm:weightedversion} Let $w_0, w_1: \R \rightarrow \R_+$ be two continuous (weight) functions satisfying
\begin{align*}
    C^{-1} \leq w_j(x) \leq C \quad \mbox{for some $C>0$ and any $x\in \R$.}
\end{align*}
Then for any $h \in \mathbb{H}^{p,q}(\mathcal{S})$ and $z=x+\ii y\in \mathcal{S}$ we have
\begin{align*}
    |h(x+\ii y)| \leq H_{p,q}(y) \exp\left(-(\mathscr{P}_y \ast \log w_0)(x) - (\mathscr{P}_{1-y} \ast \log w_1)(x)\right) \norm{h_0 w_0}_{\mL^p(\R)}^{1-y} \norm{h_1 w_1}_{\mL^q(\R)}^y, 
\end{align*}
where $\mathscr{P}_y$ is the Poisson kernel and $H_{p,q}(y)$ is given by Theorem~\ref{thm:optimum}. Moreover, this inequality is sharp.
\end{theorem}

Next, we present a natural generalization of Stein's (complex) interpolation Theorem (cf. \cite[Theorem IX.21]{RS75}). Since the proof is just the same as the proof of Stein's interpolation theorem in \cite{RS75} with Hadamard's three lines lemma replaced by Theorem~\ref{thm:optimum}, we shall skip it here.

\begin{theorem}[Generalized Stein interpolation theorem] \label{thm:rieszthorin} Let $(X,\mu)$ and $(Y,\nu)$ be $\sigma$-finite measure spaces and let $\{ T(z) \}_{z\in \mathcal{S}}$ be a family of linear operators mapping finitely simple functions in $X$ to measurable functions in $Y$. Let $1\leq p_0,q_0,p_1,q_1 \leq\infty$ and suppose that for any simple functions $f : X \rightarrow \C$ and $g: Y \rightarrow \C$, the function $z \rightarrow h_{f,g}(z) \coloneqq \inner{g,T(z) f}$ belongs to $\mathbb{H}^{r,s}(S)$ for some $1\leq r,s\leq \infty$ and satisfies the bounds
\begin{align*}
    |h_{f,g}(x)| &\leq M_0(x)\norm{g}_{L^{q_0^\ast}(Y)} \norm{f}_{L^{p_0}(X)}\\
    |h_{f,g}(x+\ii)| &\leq M_1(x)\norm{g}_{L^{q_1^\ast}(Y)} \norm{f}_{L^{p_1}(X)} , 
\end{align*}
for almost every $x \in \R$ and some measurable functions $M_0, M_1 : \R \rightarrow [0,\infty)$. Then for any $z = x+ \ii \alpha \in S$ and any $1\leq s_0,s_1 \leq \infty$, the operator $T(z)$ satisfies the bound
\begin{align*}
    \norm{T(x+\ii \alpha) f}_{L^{q_\alpha}(Y)} \leq H_{s_0,s_1}(\alpha) \norm{M_0}_{L^{s_0}(\R)}^{1-\alpha}\norm{M_1}_{L^{s_1}(\R)}^\alpha \norm{f}_{L^{p_\alpha}(X)},
\end{align*}
for any simple $f:X \rightarrow \C$, where $H_{s_0,s_1}(\alpha)$ are the values from Theorem~\ref{thm:optimum} and $p_\alpha^{-1} = (1-\alpha)p_0^{-1} + \alpha p_1^{-1}$ and $q_\alpha^{-1} = (1-\alpha)q_0^{-1} + \alpha q_1^{-1}$. In particular, if $M_0 \in L^{s_0}(\R)$ and $M_1 \in L^{s_1}(\R)$ for some $1\leq s_0,s_1\leq \infty$, then $T(x+\ii \alpha)$ extends uniquely to a bounded operator from $L^{p_\alpha}(X)$ to $L^{q_\alpha}(Y)$ satisfying the bound
\begin{align*}
    \norm{T(x+\ii \alpha)}_{L^{p_\alpha} \rightarrow L^{q_\alpha}} \leq H_{p_\alpha,q_\alpha}(\alpha) \norm{M_0}_{L^{s_0}(\R)}^{1-\alpha} \norm{M_1}_{L^{s_1}(\R)}^\alpha.
\end{align*}
\end{theorem}

\begin{remark*}[Controlled growth assumption] The assumption $h_{f,g} \in \mathbb{H}^{r,s}(S)$ in Theorem~\ref{thm:rieszthorin} can be replaced by any controlled growth assumption that guarantees both the existence of boundary values and that the following form of Cauchy's integral theorem holds:
\begin{align*}
    \int_\R h_{f,g}(x+\ii y) \phi(x+\ii y) \mathrm{d} x = \int_\R h_{f,g}(x+\ii y') \phi(x+\ii y') \mathrm{d} x \quad \mbox{for any $0<y,y'<1$,}
\end{align*}
and for any $\phi = \widehat{\psi}$ with $\psi \in C_c^\infty(\R)$. For instance, provided that the boundary values of $h_{f,g}$ exists in the sense of Definition~\ref{def:trace}, any bound of the form
\begin{align*}
    |h_{f,g}(z)|\lesssim (1+|z|)^k\biggr( \frac{1}{|\mathrm{Im}(z)|^k} + \frac{1}{|1-\mathrm{Im}(z)|^k}\biggr) \quad \mbox{for some $k \in \N$}
\end{align*}
suffices.
\end{remark*}

\begin{remark*}[Further extensions] We can similarly use Theorem~\ref{thm:optimum} (in the place of Hadamard's three lines lemma) to generalize further results concerning family of analytic operators in general interpolation spaces. For instance, one can obtain an analogous extension of the Calderón-Lions interpolation theorem \cite[Theorem IX.19]{RS75}. These applications are beyond the scope of the current paper and will not be discussed any further here.
\end{remark*}

The last applications we present are about Lieb-Thirring (LT) inequalities. These are applications of the special case $H_{\infty,2}(\alpha)$ and were already discussed in \cite{CCR24}. Here, however, we use the explicit formula for the optimal values (see Remark~\ref{rem:analyticformula}) to obtain a much simpler\footnote{As in \cite{CCR24}, the bounds from Corollary~\ref{cor:LTCLR} are implicitly stated in terms of $\mL^2$ norms of a non-trivial (and rather complicated) function.} analytic formula for the bounds derived there. 

\begin{corollary}[LT and CLR bound]\label{cor:LTCLR} Let $d \in \N$ and $0 < s < \infty$, and let $H_V = (-\Delta)^s +V$ be the generalized Schr\"odinger operator with some measurable (potential) function $V: \R^d \rightarrow \R$. Let $N_0(H_V)$ denote the number of negative eigenvalues of $H_V$ (counting multiplicity) and $\mathrm{tr} (H_V)_-$ denote minus the sum of the negative eigenvalues of $H_V$. Then the optimal constant $L_{0,d,s}$ in the CLR inequality (under the restriction $0 < s< d/2$)
\begin{align*}
    N_0\left(H_V\right) \leq L_{0,d,s} \int_{\R^d}V_-(x)^{\frac{d}{2s}} \mathrm{d} s,
\end{align*}
and the optimal constant $L_{1,d,s}$ in the LT inequality
\begin{align*}
    \mathrm{tr} \left(H_V\right)_- \leq L_{1,d,s} \int_{\R^d} V_-(x)^{\frac{2s + d}{2s}} \mathrm{d} x
\end{align*}
satisfy
\begin{align}
\frac{L_{0,d,s}}{L^{\rm cl}_{0,d}} &\leq \frac{\pi}{\alpha \sin(\pi \alpha)}\exp\left(\frac{\mathrm{CI}_2\left(2\pi(1-\alpha)\right)}{\pi \alpha}\right), \quad \mbox{with $\alpha = \frac{2s}{d}$,}
\label{eq:CLRbound}
\intertext{and}
\frac{L_{1,d,s}}{L^{\rm cl}_{1,d,s}} &\leq \frac{\pi (1-\alpha)^{\frac{1}{\alpha}}}{\alpha \sin(\pi \alpha)} \exp\left(\frac{\mathrm{CI}_2\left(2\pi(1- \alpha)\right)}{\pi \alpha}\right), \quad \mbox{with $\alpha = \frac{2s}{d+2s}$,} \label{eq:LTbound}
\end{align}
where $\mathrm{CI}_2$ is the Clausen function~\eqref{eq:clausendef} and the semi-classical constant is $L^{\rm cl}_{1,d,s} = \frac{2s}{d+2s} L^{\rm cl}_{0,d}$ with $L^{\rm cl}_{0,d} = \frac{|B_1|}{(2\pi)^d}$.
\end{corollary}

\begin{proof}
The proof follows from~\cite[Theorem 1.2]{CCR24} and~\eqref{eq:analyticformula2case}.
\end{proof}

\begin{remark*}[Best known bounds] \label{rem:LTbounds} The bound in~\eqref{eq:LTbound} leads to $L_{1,1,1}/L_{1,1,1}^{\rm cl} \leq 1.447$ (see \cite{CCR24}), which can be lifted to higher dimensions (see \cite{LW00}) and is the best known so far. For the CLR inequality, the best bound for $s=1$ in dimensions $d=3$ or $d=4$ is due to Lieb \cite{Lie76}. For all the other cases, i.e., $d\geq 5$ or $s\neq 1$ for the CLR and all $d$ with $s\neq 1$ for the LT, the bounds in Corollary~\ref{cor:LTCLR} appears to be the best available to date\footnote{For further details on Lieb-Thirring inequalities, see the reviews \cite{Fra21,Sch22,Nam21} and the book~\cite{FLW23}}. For instance, they yield the asymptotic bounds
\begin{align*}
    \limsup_{d/2s \rightarrow \infty} \frac{L_{1,d,s}}{L_{1,d,s}^{\rm cl}} \leq 4\pi^2 \ee^{-2} \approx 5.343 \quad\mbox{and}\quad \limsup_{d/2s \rightarrow \infty} \frac{L_{0,d,s}}{L_{0,d}^{\rm cl}} \leq 4\pi^2 \ee^{-3} \approx 1.966.
\end{align*}
\end{remark*}

\subsection{Outline of the paper}

We now briefly outline the key steps in the proof of our main results and how they are organized throughout the next sections. 

In Section~\ref{sec:Hardy} we collect various results on the spaces $\mathbb{H}^{p,q}(S)$. These include a characterization of the boundary values of such functions and a few technical lemmas that will be used throughout our proofs.

Section~\ref{sec:duality} contains the proof of the duality relation in Theorem~\ref{thm:duality}. Their proofs combine some tools of complex and convex analysis to first derive a dual formulation of problem~\eqref{eq:pq3lineproblem} over a suitable space of meromorphic functions, and then show, via a simple Blaschke product factorization, that this dual formulation is equivalent to the original problem with the H\"older conjugate exponents. This suffices to complete the proof of Theorem~\ref{thm:duality}.

The proof of  Theorem~\ref{thm:optimizers} is presented in Section~\ref{sec:eulerlagrange}. The first step here consists in writing down a set of coupled Euler-Lagrange equations for the pair of optimizers of problem~\eqref{eq:pq3lineproblem} and the associated dual problem over meromorphic functions described in Section~\ref{sec:duality}. These equations are non-local and non-linear in the sense that they relate the boundary values of a holomorphic function (the primal optimizer) and a meromorphic function (the dual optimizer) along distinct boundaries of the strip in a nonlinear way. Fortunately, it turns out that by combining (a weak version of) the Schwarz reflection principle and Liouville's theorem, we are able to derive an explicit formula for the product between the primal and dual optimizers (Lemma~\ref{lem:dualproduct}). This formula depends only on $\alpha \in (0,1)$; moreover, it justifies the Ansatz $h= \ee^{\phi}$ for the optimizers and allows us to construct the function $\phi$ via the well-known Poisson integral representation for harmonic functions on the strip. This considerably simplifies the construction for the optimizers in \cite{CCR24}, which relied on the Fourier transform of singular value distributions and was very specific for the case $p=\infty$ and $q=2$. 

Some additional results that are used throughout the proofs are collected in the appendix. More precisely, in Appendix~\ref{app:analyticmeasures} we prove that measures whose Fourier transform has an one-sided exponential decay in an averaged sense are absolutely continuous with respect to the Lebesgue measure. This result is used in Section~\ref{sec:Hardy} to prove the characterization of the boundary values of functions in $\mathbb{H}^{p,q}(S)$ when either $p=1$ or $q=1$. In Appendix~\ref{app:Poisson} we collect some well-known facts about the Poisson kernel on the strip that are used throughout the paper. In Appendix~\ref{app:analyticformula}, we present a proof of formula~\eqref{eq:analyticformula}. Finally, in Appendix~\ref{app:weigthed} we prove Theorem~\ref{thm:weightedversion}.
\addtocontents{toc}{\protect\setcounter{tocdepth}{1}}
\section{Hardy-like spaces on the strip} \label{sec:Hardy}

In this section, we introduce the mixed Hardy spaces $\mathbb{H}^{p,q}(S)$ and derive some of its basic properties. Some of the results presented here are extensions of the results obtained in \cite{CCR24} for the special case $\mathbb{H}^{p,2}(S)$. Therefore, we shall sporadically refer to that paper for some details.

Let us start with some notation. First, we recall that $S$ denotes the infinite horizontal open strip
\begin{align*}
    S = \{z = x+iy \in \C : 0 < y < 1 \} \subset \C  .
\end{align*}
More generally, we denote by $S_{a,b}$ the strip
 \begin{align*}
     S_{a,b} \coloneqq \{ z = x + \ii y \in \C: a<y<b\} \quad \mbox{for some $-\infty \leq a < b \leq \infty$.}
 \end{align*}
 In addition, we use the notation $S_{a,b}^\ast$ for the reflected (via complex conjugation) strip 
 \begin{align*}
     S_{a,b}^\ast \coloneqq \{\overline{z} : z \in S_{a,b} \} = \{z+\ii y : -b<y<-a\} = S_{-b,-a}.
 \end{align*}
 The space of all holomorphic functions on $S_{a,b}$ is denoted by $\mathbb{H}(S_{a,b})$. Moreover, for any $h\in \mathbb{H}(S_{a,b})$, we denote by $h^\ast\in\mathbb{H}(S_{a,b}^\ast)$ the reflected function
\begin{align}
    h^\ast(z) \coloneqq \overline{h(\overline{z})}, \quad z \in S_{a,b}^\ast. \label{eq:reflecteddef}
\end{align}
Let us also recall the notation
\begin{align*}
	h_y(x) \coloneqq h(x+\ii y)\quad \mbox{for $a<y<b$ and $h \in \mathbb{H}(S_{a,b})$.}
\end{align*} 
We can now define the mixed Hardy spaces $\mathbb{H}^{p,q}(S_{a,b})$ in a general strip $S_{a,b}$ as follows:
\begin{definition}[$\mathbb{H}^{p,q}(S_{a,b})$ spaces] \label{def:Hardy}
Let $1\leq p \leq \infty$, we denote by $\mathbb{H}^{p,q}(S_{a,b})$ the space of holomorphic functions $h:S \rightarrow \C$ satisfying the bound
\begin{align}
\norm{h}_{\mathbb{H}^{p,q}} \coloneqq \sup_{a<y<b} \inf_{f+g = h_y} \{(b-y)^{-1} \norm{f}_{L^p(\R)} + (y-a)^{-1} \norm{g}_{L^q(\R)}\} <\infty . \label{eq:holomorphicbound}
\end{align}
Moreover, we say that $h \in \mathbb{H}^{p,q}(S_{a,b}^\ast)$ if the reflected function $h^\ast$ (see~\eqref{eq:reflecteddef}) belongs to $\mathbb{H}^{p,q}(S_{a,b})$.
\end{definition}

\begin{remark*}[Classical Hardy spaces on the strip] Note that for $p=q$,
\begin{align*}
\sup_{0<y<1} \norm{h_y}_{L^p(\R)} \leq \norm{h}_{\mathbb{H}^{p,p}(S)}  \leq 2 \sup_{0<y<1} \norm{h_y}_{L^p(\R)}.
\end{align*}
Hence $\mathbb{H}^{p,p}(S)$ is the standard (translation invariant) $p$-Hardy space on the strip (cf. \cite{BK07}).
\end{remark*}

Next, let us introduce the definition of tempered boundary values of a holomorphic function on the strip. For this, we use the convention
\begin{align}
    \inner{f,g} = \int_\R f(x) \overline{g(x)} \mathrm{d} x \label{eq:innerproductdef}
\end{align}
for the $L^2(\R)$ inner-product of two measurable functions $f,g:\R \rightarrow \C$.
\begin{definition}[Tempered trace] \label{def:trace} Let $h : S_{a,b} \rightarrow \C$ be a holomorphic function on the strip $S_{a,b}$ for some $-\infty<a < b < \infty$, then we say that $h$ admits a trace (or has boundary values) along the line $
R+a\ii$, if there exists a tempered distribution $T \in \mathcal{S}^\prime(\R)$ such that
\begin{align}
    \lim_{y \downarrow a} \inner{h_y,\phi} = T(\overline{\phi}),  \quad \mbox{for any Schwartz function $\phi \in \mathcal{S}(\R)$,} \label{eq:boundaryvalues}
\end{align}
 where $\overline{\phi}$ denotes the complex-conjugated function\footnote{The seemingly unnatural complex-conjugation here is only necessary to match the convention adopted for $\inner{\cdot,\cdot}$ in \eqref{eq:innerproductdef} and justify the identification $T =h_a$ whenever the boundary trace is a function.} $\overline{\phi}(x) = \overline{\phi(x)}$. Similarly, we say that $h$ admits a trace $S\in \mathcal{S}^\prime(\R)$ along the line $\R+\ii b$ if
\begin{align*}
    \lim_{y \uparrow b} \inner{h_y,\phi} = S(\overline{\phi}),  \quad \mbox{for any Schwartz function $\phi \in \mathcal{S}(\R)$.} 
\end{align*}
In this case, we denote the distributions $T$ and $S$ respectively by $h_a$ and $h_b$.
\end{definition}

A key fact about the space $\mathbb{H}^{p,q}(S)$ is that its boundary values along the real-axis are well-defined and can be identified with the subspace of functions in $L^p(\R)$ whose Fourier transform has an one-sided exponential decay in a suitable sense. To state this result precisely, we introduce the following additional notation.

For $\mu \in \R$, we denote by $\exp_\mu$ the exponential function $k \in \R \mapsto \exp_\mu(k) = \ee^{\mu k}$. (If $\mu=1$ we write simply $
\exp$.) Then for any tempered distribution $T \in \mathcal{S}^\prime(\R)$, we denote by $T \exp_\mu \in \mathcal{D}^\prime(\R)$  the product of $T$ with the exponential function $\exp_\mu$, i.e., 
\begin{align*}
	\bigr(T \exp_\mu\bigr)(\phi) = T(\exp_\mu \phi), \quad \mbox{for any $\phi \in C_c^\infty(\R)$.}
\end{align*}
If $T \exp_\mu$ extends to a tempered distribution, we define its inverse Fourier transform by
\begin{align}
	\mathcal{F}^{-1}(T \exp_{\mu})(\phi) = T(\exp_\mu \mathcal{F}^{-1}\phi), \,\,\,\,\mbox{with the convention}\,\,\,\, \mathcal{F}^{-1}\phi(k) = \frac{1}{2\pi} \int_\R \phi(x) \ee^{\ii xk} \mathrm{d}x  \label{eq:FTconvention}
\end{align} 
for any function $\phi$ in the Schwartz space $\mathcal{S}(\R)$. Note that since $C_c^\infty(\R)$ is dense in the Schwartz space, the extension above (if existing) is unique.

We can now state the characterization for the boundary values of $\mathbb{H}^{p,q}(S)$ as follows.
\begin{theorem}[Boundary values of $\mathbb{H}^{p,q}(S)$]\label{thm:holomorphicextension} For any $v\in \{f \in L^p(\R) : \mathcal{F}^{-1}(\widehat{f}\exp_{-1} )\in L^q(\R)\}$, there exists an unique function $h\in \mathbb{H}^{p,q}(S)$ such that $h_0 = v$. In this case, we have $\widehat{h}_1 = \widehat{v} \exp_{-1}$. Conversely, for any $h \in \mathbb{H}^{p,q}(S)$ we have $h_0 \in \{f \in L^p(\R) : \mathcal{F}^{-1}(\widehat{f} \exp_{-1}) \in L^q(\R)\}$. Moreover, we have the equivalence of norms
\begin{align}
    \norm{h}_{\mathbb{H}^{p,q}(S)} \lesssim \norm{h_0}_{L^p(\R)} + \norm{h_1}_{L^q(\R)} \lesssim \norm{h}_{\mathbb{H}^{p,q}(S)}. \label{eq:normequivalence}
\end{align}
\end{theorem}

To prove Theorem~\ref{thm:holomorphicextension}, we shall need a few preliminary lemmas. The first one is a pointwise bound on $h$ that is uniform away from the boundary. The proof of this result is a simple adaptation of a standard mean value argument used in the theory of Hardy spaces (see, e.g.,\cite[Lemma 11.3]{Mas09}). For the proof details, we refer the reader to \cite[Lemma B.1]{CCR24}.  

\begin{lemma}[Uniform estimate]\label{lem:uniformest} Let $1\leq p,q \leq \infty$ and $h \in \mathbb{H}^{p,q}(S)$, then we have
\begin{align*}
|h(z)| \lesssim  \norm{h}_{\mathbb{H}^{p,q}} \bigr(y^{-\frac{1}{p}} + (1-y)^{-\frac{1}{q}}\bigr),\quad \mbox{for any $0<y<1$,}
\end{align*}
where $\norm{h}_{\mathbb{H}^{p,q}}$ is defined in \eqref{eq:holomorphicbound}.
\end{lemma}


As an immediate consequence of the above lemma, we have the following result.
\begin{proposition}[Banach space] For any $1\leq p,q\leq \infty$, the space $\mathbb{H}^{p,q}(S)$ endowed with the norm $\norm{\cdot}_{\mathbb{H}^{p,q}}$ defined in \eqref{eq:holomorphicbound} is a Banach space.
\end{proposition}
\begin{proof}
    It is immediate to verify that $\norm{\cdot}_{\mathbb{H}^{p,q}}$ is a norm. To show that $\mathbb{H}^{p,q}$ is complete, we observe that, by Lemma~\ref{lem:uniformest}, any Cauchy sequence $\{h^{(n)}\}_{n\in \N} \subset \mathbb{H}^{p,q}(S)$ is also a Cauchy sequence with respect to uniform convergence in compact subsets of $S$. In particular, it must converge locally uniformly to a holomorphic function $h \in \mathbb{H}(S)$. That $h \in \mathbb{H}^{p,q}(S)$ and $\norm{h-h^{(n)}}_{\mathbb{H}^{p,q}} \rightarrow 0$ follows from a standard triangle inequality argument. 
\end{proof}

We now prove the existence of boundary values (in the weak sense) for $\mathbb{H}^{p,q}$ functions. 

\begin{lemma}[Existence of boundary values in $\mathbb{H}^{p,q}(S)$]\label{lem:existenceboundary} Let $h\in \mathbb{H}^{p,q}(S)$, then $h$ admits boundary values along the lines $\R$ and $\R+\ii$ in the sense of Definition~\ref{def:trace}. Moreover, the boundary values satisfy 
\begin{align}
    \norm{h_0}_{L^p(\R)} + \norm{h_1}_{L^q(\R)} \leq 2 \norm{h}_{\mathbb{H}^{p,q}(S)}\label{eq:firstnormbound}
\end{align}
and $\widehat{h}_1 =\exp_{-1} \widehat{h}_0$.
\end{lemma}

\begin{proof} The proof here follows the same steps in the proof of \cite[Lemma 3.1]{CCR24} with only minor modifications. Nevertheless, we shall present the full proof because similar arguments will be used later.

The first step is to show that for any $\phi \in C_c^\infty(\R)$, the function
\begin{align}
 \inner{h_y, \widehat{\phi}_{y_0-y}} = c(y_0)\quad \mbox{ is independent of $y \in (0,1)$.} \label{eq:closedcurveintegral}
\end{align}
For this, we note that the function $z\mapsto h^\ast(z) \widehat{\phi}(z + iy_0)$ is holomorphic on $S^\ast$ (as $\widehat{\phi}$ is entire), and therefore,
\begin{align*}
	\int_{-R}^R \biggr(\overline{h_y(x)} \widehat{\phi}_{y_0-y}(x)  - \overline{h_{y'}(x)} \widehat{\phi}_{y_0-y'}(x)\biggr)\mathrm{d} x  + \int_{y'}^y \biggr(\overline{h_{w}(-R)}\widehat{\phi}_{y_0-w}(-R) - \overline{h_{w}(R)}\widehat{\phi}_{y_0-w}(R)\biggr) \mathrm{d} w = 0
\end{align*}
by Cauchy's integral theorem. In the limit $R\ra \infty$, the first term converges to $\inner{h_y, \widehat{\phi}_{y_0-y}} - \inner{h_{y'},\widehat{\phi}_{y_0-y'}}$ by dominated convergence. On the other hand, the second term vanishes thanks to the uniform (with respect to $R$) control on $h_w(R)$ from Lemma~\ref{lem:uniformest} and the fast decay of $\phi_w(R)$. This completes the proof of~\eqref{eq:closedcurveintegral}.

The next step is to show the existence of a unique weak limit for $h_y$ as $y$ approaches the boundary of the strip. For this, we let $f+g = h$ be a decomposition satisfying
\begin{align}
\norm{f_y}_{L^p(\R)} \lesssim (1-y)\quad\mbox{and}\quad \norm{g_y}_{L^q(\R)} \lesssim y, \label{eq:fgdecomposition}
\end{align} 
which exists since $\norm{h}_{\mathbb{H}^{p,q}} < \infty$. (Note that we do not require $f$ and $g$ to be holomorphic.) Then from Banach-Alaoglu, we can extract subsequences $y_n \ra 0$, $v_n \ra 1$  and $h_0 \in L^p(\R)$ (or $\mathcal{M}(\R)$ if $p=1$) and $h_1 \in L^q(\R)$ (or $\mathcal{M}(\R)$ if $q=1$) such that $f_{y_n} \ra h_0$ and $g_{v_n} \ra h_1$ in the respective weak star topologies. We can now prove that the limit is unique by using \eqref{eq:closedcurveintegral}. Precisely, from \eqref{eq:closedcurveintegral}
\begin{align}
	\inner{h_y,\widehat{\phi}_0} - \inner{h_0, \widehat{\phi}_0} &= \inner{h_y,\widehat{\phi}_{-y}} - \inner{h_0, \widehat{\phi}_0}  - \inner{h_y, \widehat{\phi}_{-y}-\widehat{\phi}_0} \nonumber \\
	&=  \inner{h_{y_n}, \widehat{\phi}_{-y_n}} - \inner{h_0,\widehat{\phi}_0} + \inner{h_y, \widehat{\phi}_{-y}-\widehat{\phi}_0} \nonumber\\
	&= \inner{h_{y_n}-h_0, \widehat{\phi}_0} + \inner{h_{y_n}, \widehat{\phi}_{-y_n}-\widehat{\phi}_0} - \inner{h_y,\widehat{\phi}_{-y}-\widehat{\phi}_0}  \label{eq:auxest}
\end{align}
for any $\phi \in C^\infty_c(\R)$. But from \eqref{eq:fgdecomposition}, we have $\inner{h_{y_n},\widehat{\phi}_{0}} = \inner{f_{y_n} + g_{y_n},\widehat{\phi}_0} \ra \inner{h_0,\widehat{\phi}_0}$ as $n\ra \infty$; therefore, the first term can be made arbitrarily small by choosing $n$ large. On the other hand, since 
\begin{align*}
    \sup_{0<y<1} \norm{h_y}_{L^p+L^q} \lesssim \norm{h}_{\mathbb{H}^{p,q}}< \infty \quad \mbox{ and }\quad \widehat{\phi}_{-y} \ra \widehat{\phi}_0\quad \mbox{strongly in $L^1 \cap L^\infty$ as $y\ra 0$,}
\end{align*}
the second term in \eqref{eq:auxest} can also be made arbitrarily small by choosing $n$ large and the third term vanishes in the limit $y\ra 0$. For the case $y\ra 1$, the exact same argument holds. In summary, \eqref{eq:auxest} holds for Fourier transform of compactly supported smooth functions. As the latter are dense in the Schwartz space, eq.~\eqref{eq:boundaryvalues} holds for any Schwartz function.

The proof for the case $p,q>1$  is essentially complete at this point. For the case where $p=1$ and/or $q=1$, we need to show that the limit measures $h_0 \in \mathcal{M}(\R)$ and/or $h_1 \in \mathcal{M}(\R)$ are absolutely continuous with respect to the Lebesgue measure. For this, we note that by the preceding arguments, \eqref{eq:closedcurveintegral} extends to $y=0$ and $y=1$, i.e.,
\begin{align}
	\inner{h_0, \widehat{\phi}_0} = \inner{h_1, \widehat{\phi}_{-1}} . \label{eq:extendedlineintegral}
\end{align}
Thus, since $\widehat{\phi}_{-1}$ is the Fourier transform of $\phi  \exp_{-1}$, i.e.,
\begin{align*}
    \widehat{\phi}_{-1}(k) = \widehat{\phi}(k-\ii) = \int_\R \phi(x) \ee^{-\ii x(k-\ii)} \mathrm{d} x = \int_\R \phi(x) \exp_{-1}(x) \ee^{-\ii k x} \mathrm{d} x = \widehat{\phi \exp_{-1}}(x),
\end{align*}
from \eqref{eq:extendedlineintegral} and recalling the convention adopted for the inverse F.T. in \eqref{eq:FTconvention}, we have
\begin{align*}
    \left(\mathcal{F}^{-1}h_0\right) \left(\overline{\phi}\right) &= h_0\left(\mathcal{F}^{-1}\left(\overline{\phi}\right)\right) = \frac{1}{2\pi} \inner{h_0, 
 \widehat{\phi_0}} = \frac{1}{2\pi} \inner{h_1, \widehat{\phi}_{-1}} =  h_1\left(\mathcal{F}^{-1}\left(\overline{\phi} \exp_{-1}\right)\right)
\end{align*}
This implies that $\mathcal{F}^{-1}(h_0) = \mathcal{F}^{-1}(h_1) \exp_{-1}$, or equivalently, $\widehat{h}_1 = \exp_{-1} \widehat{h}_0$. As both $h_0$ and $h_1$ belongs to $\mathcal{M}(\R) + L^\infty(\R)$ (independent of $p,q$), the absolute continuity of $h_0$ and $h_1$ follows from Lemma~\ref{lem:abscont}.

Inequality~\eqref{eq:firstnormbound} follows from the lower semi-continuity of the norms $\norm{\cdot}_{L^p(\R)}$ with respect to weak or weak$^\ast$ convergence.
\end{proof}

We can now finish the proof of Theorem~\ref{thm:holomorphicextension}.

\begin{proof}[Proof of Theorem~\ref{thm:holomorphicextension}] Let $v \in L^p(\R)$ with $w \coloneqq \mathcal{F}^{-1}(\widehat{v} \exp_{-1}) \in L^q(\R)$, then the idea is to define 
\begin{align}
 &h(z) = \underbrace{\left(v\ast \mathscr{P}_y\right)(x)}_{\coloneqq f_y(x)} + \underbrace{\left(\mathscr{Q}_y \ast {w}\right)(x)}_{\coloneqq g_y(x)}  \label{eq:hdef}
 \end{align}
 where $\mathscr{P}_y$ and $\mathscr{Q}_y$ are approximate identities when $y\downarrow 0$ and $y\uparrow 1$, respectively, and they satisfy
 \begin{align}
&\ee^{\ii xk} (\widehat{\mathscr{P}}_y(k) + \widehat{\mathscr{Q}}_y(k)\ee^{-k}) =  \ee^{\ii zk} \label{eq:condition1}
\intertext{and}
&\norm{\mathscr{P}_y}_{L^1(\R)} \lesssim 1-y \quad\mbox{and}\quad \norm{\mathscr{Q}_y}_{L^1(\R)} \lesssim y. \label{eq:conditions}
\end{align}
Indeed, if estimate~\eqref{eq:conditions} holds, then we can apply Young's convolution inequality to obtain 
\begin{align}
    \norm{h}_{\mathbb{H}^{p,q}} \leq |y|^{-1}\norm{f_y}_{L^p} +(1-y)^{-1} \norm{g_y}_{L^q} \lesssim \norm{v}_{L^p} + \norm{w}_{L^q}. \label{eq:normest2}
\end{align}
Hence, by integrating against Schwartz functions $\phi$ and using the approximate identity property of $\mathscr{P}_y$ and $\mathscr{Q}_y$ we find
\begin{align*}
	\inner{h_y, \phi} &=  \inner{\mathscr{P}_y \ast v, \phi} + \inner{\mathscr{Q}_y \ast {w}, \phi} \ra \begin{dcases} \inner{v, \phi} = \inner{v, \phi}, \quad \mbox{as $y \downarrow 0$,} \\
	\inner{w, \phi} = \inner{w,\phi},\quad \mbox{as $y \uparrow 1$.} \end{dcases}
\end{align*}
which implies $h_0 = v$ and $h_1 = {w}$. 

The choice of $\mathscr{P}_y$ and $\mathscr{Q}_y$ is to some extent arbitrary. For concreteness (and since these functions will naturally appear later), we choose
\begin{align*}
\widehat{\mathscr{P}}_y(k) = \frac{\sinh\left((1-y)k\right)}{\sinh(k)},\quad \mbox{and}\quad \widehat{\mathscr{Q}}_y = \widehat{\mathscr{P}}_{1-y}(k) = \frac{\sinh(y k)}{\sinh(k)}.
\end{align*}
Then by Lemma~\ref{lem:PoissonFourier}, we have
\begin{align*}
    \mathscr{P}_y(x) =\frac{1}{2} \frac{\sin(\pi y)}{\cosh(\pi x) - \cos(\pi y)} \quad \mbox{and} \quad \mathscr{Q}_y(x) = \frac{1}{2} \frac{\sin(\pi y)}{\cosh(\pi x) + \cos(\pi y)}, 
\end{align*}
In particular, $\mathscr{P}_y$ is an approximate identity as $y\downarrow 0$ by Lemma~\ref{lem:PoissonRepresentation},  and estimates~\eqref{eq:conditions} holds by Lemma~\ref{lem:PoissonFourier}. It remains to show that $h(z)$ is holomorphic. 

For this, we let $\widehat{\phi}$ be a smooth compactly supported function satisfying $\widehat{\phi}(0) = 1$ and set 
\begin{align*}
    \widehat{v}^\epsilon \coloneqq \widehat{v} \widehat{\phi}(\epsilon \cdot)\quad \mbox{and}\quad \widehat{w}^\epsilon \coloneqq \widehat{v} \exp_{-1} \widehat{\phi}(\epsilon \cdot), \quad \mbox{for $\epsilon >0$.} 
\end{align*}
Then $\widehat{v}^\epsilon$ is a tempered distribution with compact support. Consequently, $v^\epsilon(z)$ is well-defined and entire by the Paley-Wiener theorem. Moreover, by~\eqref{eq:condition1} it satisfies
\begin{align}
    v^\epsilon(z) &= \frac{1}{2\pi} \int_\R \ee^{\ii z k} \widehat{v}(k) \widehat{\phi}(\epsilon k) \mathrm{d} k = \frac{1}{2\pi} \int_\R \ee^{\ii x k} \widehat{\mathscr{P}}_y(k) \widehat{v}^\epsilon( k)\mathrm{d} k + \frac{1}{2\pi} \int_\R \ee^{\ii x k} \widehat{\mathscr{Q}}_y(k) \widehat{w}^\epsilon(k) \mathrm{d} k\nonumber\\
    &= \left(\mathscr{P}_y \ast v^\epsilon\right)(x) + \left(\mathscr{Q}_y \ast w^\epsilon\right)(x), \quad \mbox{for any $z= x+\ii y \in S$.} \label{eq:funid}
\end{align}
Thus from Young's convolution inequality, the sequence $\{v^\epsilon\}_{\epsilon >0}$ is uniformly bounded in $\mathbb{H}^{p,q}(S)$. In particular, by Lemma~\ref{lem:uniformest}, this sequence is also uniformly (pointwise) bounded on compact subsets of $S$. Cauchy's integral formula then implies that the set $\{v_\epsilon\}_{\epsilon>0}$ is also equicontinuous in compact subsets of $S$. Thus by Arzela-Ascoli's theorem\footnote{This application of Arzela's Ascoli theorem is also known as Montel's theorem in complex analysis.}, we can extract a subsequence that converges locally uniformly to some holomorphic function on the strip, which we denote by $\widetilde{h}$. 

To conclude, we now show that $\widetilde{h}_y(x) = h_y(x)$ for almost every $x\in \R$,. For this, we note that from the approximate identity property of $\phi^\epsilon(x) = \phi(x/\epsilon)/\epsilon$ and Young's convolution inequality we have
\begin{align*}
    \lim_{\epsilon \downarrow 0} \mathscr{P}_y \ast v^\epsilon = \lim_{\epsilon \downarrow 0} (\mathscr{P}_y \ast \phi^\epsilon) \ast v = \mathscr{P}_y \ast v \quad \mbox{and}\quad \lim_{\epsilon \downarrow0} \mathscr{Q}_y\ast w^\epsilon = \lim_{\epsilon \downarrow 0} (\mathscr{Q}_y \ast \phi^\epsilon) \ast w = \mathscr{Q}_y \ast w, 
\end{align*}
where the convergence is strong in $L^p(\R)$ and $L^q(\R)$, respectively. Therefore, the result follows from~\eqref{eq:funid}, the definition of $h$ in~\eqref{eq:hdef}, and the fact that $v^\epsilon(z) \ra \widetilde{h}(z)$ for any $z\in S$.

The converse statement in Theorem~\ref{thm:holomorphicextension} was proved in Lemma~\ref{lem:existenceboundary}, and the equivalence of norms in~\eqref{eq:normequivalence} follows from estimates~\eqref{eq:firstnormbound} and \eqref{eq:normest2}.  \end{proof}

\begin{remark*}[Poisson representation of holomorphic functions on the strip] The proof of Theorem~\ref{thm:holomorphicextension} shows that the Poisson representation $u(x+ \ii y) = \mathscr{P}_y \ast v(x) + \mathscr{P}_{1-y} \ast w(x)$, which corresponds to the solution\footnote{which is unique only among a class of harmonic functions with restricted growth, see \cite{Wid61}.} of Laplace's equation in the strip with boundary values $v$ and $w$ (see Lemma~\ref{lem:PoissonRepresentation}), is holomorphic if and only if $\widehat{v} = \widehat{w} \exp$. This result is the analogue \emph{on the strip} for the Poisson representation of Hardy functions on the half-plane (cf. \cite[Corollary 3.2]{Gar06} and \cite[Theorem 13.2]{Mas09}). 
\end{remark*}

\begin{remark*}[Strong, weak, and tempered boundary values] The proof of Lemma~\ref{lem:existenceboundary} shows that for any $h \in \mathbb{H}^{p,q}(S)$, the  convergence of $h_y$ to its boundary values takes place in the weak topology of $L^p(\R)$ (or weak-$\star$ topology for $L^\infty(\R)$) which is stronger than the convergence in the tempered distribution sense\footnote{The existence of boundary values in the tempered distribution sense can be proved (see, e.g., \cite{Til61,CKP07}) for the larger class of function satisfying the slow growth condition
\begin{align*}
 |h(z)| \lesssim |z|^m (|y|^{-n}+|1-y|^k) \quad\mbox{for some $m,n,k\geq 0$ and $y$ close to $0$.}
\end{align*}}  defined in \eqref{eq:boundaryvalues}. In fact, if $h_0 \in L^p(\R)$ and $h_1 \in L^q(\R)$ for some $q\leq p<\infty$, then the convergence to $h_0$ and $h_1$ holds in the strong $L^p(\R)$ and in the strong $L^q(\R) + L^\infty(\R)$ topologies, respectively.
\end{remark*}

Let us end this section with a few simple lemmas that will be useful later.



\begin{lemma}[Necessary and sufficient conditions for holomorphic extension]\label{lem:boundarycriteria} Let $1\leq p,q\leq \infty$ and $f\in L^{p^\ast}(\R)$, where $p^\ast$ is the H\"older conjugate of $p$. Then there exists a (unique) $F \in \mathbb{H}^{p^\ast,q^\ast}(S^\ast)$ with $F_0 =f$ if and only if there exists a constant $C>0$ such that
\begin{align}
 |\inner{f, h_0}| \leq C \norm{h_1}_{L^q(
 R)} \quad \mbox{for any $h \in \mathbb{H}^{p,q}(S)$.} \label{eq:boundarycriteria}
\end{align}
Moreover, in this case the optimal constant in the inequality above is $C= \norm{F_{-1}}_{L^{q^\ast}}$.
\end{lemma}

\begin{proof} If $f = F_0$ for some $F \in \mathbb{H}^{p^\ast,q^\ast}(S^\ast)$, then from \eqref{eq:extendedlineintegral} we have
\begin{align*}
\inner{f, h_0} = \inner{F_0, h_0} = \inner{F_{-1}, h_1} \leq \norm{F_{-1}}_{L^{q^\ast}} \norm{h_1}_{L^q},
\end{align*}
which proves the first direction. For the other direction, note that since $f\in L^p(\R)$, its Fourier transform is well-defined in the tempered distribution sense. Thus for any $\phi \in C_c^\infty(\R)$, by Plancherel's formula we have
\begin{align*}
	\inner{\widehat{f} \exp, \phi} = \inner{\widehat{f}, \phi \exp } =  2\pi \inner{f,\mathcal{F}^{-1}(\phi \exp)} = 2\pi \inner{f, (\mathcal{F}^{-1}\phi)_{-1}} \stackrel{\eqref{eq:boundarycriteria}}{\lesssim} \norm{\widehat{\phi}_0}_{L^q},
\end{align*}
where $(\mathcal{F}^{-1} \phi)_y(x) = \mathcal{F}^{-1}\phi(x+\ii y)$ is well-defined since $\mathcal{F}^{-1}\phi$ is analytic. In particular, if $1 \leq q < \infty$, we conclude that the inverse F.T. of $\widehat{f} \exp$ belongs to $L^{q^\ast}(\R)$. Theorem~\ref{thm:holomorphicextension} then completes the proof. For the case $q= \infty$, note that $\widehat{C}_c^\infty(\R)$ is dense in $C_0(\R)$ so by the Riesz representation theorem, we have $\mathcal{F}^{-1}(\widehat{f} \exp ) \in \mathcal{M}(\R)$. By Lemma~\ref{lem:abscont}, this measure must be absolutely continuous and the proof follows from Theorem~\ref{thm:holomorphicextension}.
\end{proof}

\begin{lemma}[Multiplication property]\label{lem:product} Let $h \in \mathbb{H}^{p,q}(S)$ and $m \in \mathbb{H}^{p^\ast,q^\ast}(S)$, then the product $F(z) \coloneqq m(z) h(z)$ belongs to the Hardy space $\mathbb{H}^1(S)$.
\end{lemma}

\begin{proof} For $\epsilon >0$ define $\eta_\epsilon(x) \coloneqq \eta(\epsilon x)$ for some mollifier $\eta \in C_c^\infty(\R)$ satisfying $\eta(0)=1$ and $\widehat{\eta} \geq 0$. Then the functions 
\begin{align}
    h^\epsilon = \mathcal{F}^{-1}(\widehat{h}_0 \eta^\epsilon)\quad \mbox{and}\quad m^\epsilon \coloneqq \mathcal{F}^{-1}(\widehat{m}_0 \eta^\epsilon)
\end{align}
are entire by the Paley-Wiener theorem and satisfy  
\begin{align}
    \sup_{\epsilon >0} \norm{h^\epsilon}_{\mathbb{H}^{p,q}(S)} < \infty \quad \mbox{and} \quad \sup_{\epsilon >0} \norm{m^\epsilon}_{ \mathbb{H}^{p^\ast,q^\ast}(S)} < \infty \label{eq:suphardybound} \end{align} 
by Young's convolution inequality. 

Now let $F^\epsilon$ be the entire function defined as $F^\epsilon(z) \coloneqq h^\epsilon(z) m^\epsilon(z)$. Thanks to \eqref{eq:suphardybound} and Lemma~\ref{lem:uniformest}, we have uniform control on $F^\epsilon  = m^\epsilon h^\epsilon$ away from the boundary of $S$. Therefore, we can use the same arguments from the proof of Lemma~\ref{lem:existenceboundary} (see \eqref{eq:closedcurveintegral}) to conclude that
\begin{align}
    \inner{m^\epsilon_0 h^\epsilon_0, \widehat{\phi}_1} = \inner{F^\epsilon_0,\widehat{\phi}_1} = \inner{F^\epsilon_y,\widehat{\phi}_{1-y}} = \inner{m^\epsilon_y h^\epsilon_y, \widehat{\phi}_{1-y}} \quad \mbox{for any $0\leq y \leq 1$,} \label{eq:niceequation}
\end{align}
for any $\phi \in C^\infty_c(\R)$. Hence, setting $y=1$ in~\eqref{eq:niceequation}, taking the limit $\epsilon \ra 0$, and using H\"older's inequality, we find that
\begin{align*}
    \inner{h_0 m_0 , \widehat{\phi}_1} = \inner{h_1 m_1 ,\widehat{\phi}_0} \leq \norm{h_1}_{L^q}\norm{m_1}_{L^{q^\ast}} \norm{\widehat{\phi}_0}_{L^\infty} \quad \mbox{for any $\phi \in C_c^\infty(\R)$.}
\end{align*}
Therefore, by Lemma~\ref{lem:boundarycriteria}, there exists $\widetilde{F} \in \mathbb{H}^{1,1}(S) = \mathbb{H}^1(S)$ satisfying $\widetilde{F}_0 = m_0 h_0$ and $\widetilde{F}_1 = m_1 h_1$. As~\eqref{eq:niceequation} holds for any $0\leq y\leq 1$, we conclude that $\widetilde{F}(z) = F(z)$ for any $z \in S$, which completes the proof.
\end{proof}

\begin{lemma}[Schwarz reflection principle on Hardy space]\label{lem:schwarzreflection} Let $h \in \mathbb{H}^{p,q}(S)$ be a holomorphic function satisfying 
\begin{align*}
    h_0(x) \in \R \quad \mbox{for almost every $x\in \R$.}
\end{align*}
Then $h_0$ is analytic in $\R$ and $h$ can be holomorphically extended to a function $h \in \mathbb{H}^{q,q}(S_{-1,1})$ satisfying
\begin{align*}
    h^\ast(z) = h(z) \quad \mbox{for $z \in S_{-1,1}$.} 
\end{align*}
\end{lemma}

\begin{proof} Let $h^\ast(z) = \overline{h(\overline{z})}$ for $z\in S^\ast$. Then, since $h \in \mathbb{H}^{p,q}(S)$, we have $h^\ast \in \mathbb{H}^{p,q}(S^\ast)$. Moreover, from~\eqref{eq:closedcurveintegral}, we have
\begin{align*}
    \inner{h^\ast_{-1},\widehat{\phi}_{0}} = \inner{h^\ast_{0},\widehat{\phi}_{-1}}  \quad \mbox{for any $\phi \in C_c^\infty(\R)$.}
\end{align*}
Since $h_0(x)$ is real-valued for almost every $x\in \R$, we have $h^\ast_0 = h_0$. Thus by~\eqref{eq:closedcurveintegral} again, we have
\begin{align*}
    \inner{h^\ast_{-1}, \widehat{\phi}_0} = \inner{h^\ast_0, \widehat{\phi}_{-1}} = \inner{h_0,\widehat{\phi}_{-1}} = \inner{h_1,\widehat{\phi}_{-2}} \lesssim \norm{\widehat{\phi}_{-2}}_{L^{q^\ast}(\R)},
\end{align*}
for any $\phi \in C_c^\infty(\R)$ where $1/q^\ast + 1/q =1$. The result now follows from Lemma~\ref{lem:boundarycriteria} applied to the larger strip $S_{-1,1}$ instead of $S$.
\end{proof}
\section{Duality}
\label{sec:duality}

In this section, we present a dual formulation of our variational problem over a suitable space of meromorphic functions. We then use this dual formulation to prove the natural duality relations from Theorem~\ref{thm:duality}. 

\subsection{A dual formulation over meromorphic functions}
In order to obtain a dual formulation of the variational problem~\eqref{eq:pq3lineproblem}, we first note that for any $h \in \mathbb{H}^{p,q}(S)$, the trial state
\begin{align}
	(\tau_{\tau,\beta,\omega}h)(z) = \beta e^{i \omega z} h(z-\tau), \quad \tau,\omega \in \R, \beta \in \C\setminus \{0\}, \label{eq:symmetry}
\end{align}
satisfies
\begin{align*}
	\frac{\norm{\tau_{\tau,\beta,\omega} h_{\alpha}}_{L^\infty}}{\norm{\tau_{\tau,\beta,\omega}h_0}_{L^p}^{1-\alpha} \norm{\tau_{\tau,\beta,\omega} h_1}^\alpha_{L^q}} = \frac{\norm{h_\alpha}_{L^\infty}}{\norm{h_0}_{L^p}^{1-\alpha} \norm{h_1}_{L^q}^\alpha}.
\end{align*}
In other words, problem~\eqref{eq:pq3lineproblem} is invariant under the three-parameter family of transformations in \eqref{eq:symmetry}.
These symmetries allow us to re-state problem~\eqref{eq:pq3lineproblem} as a strictly convex optimization problem and apply the celebrated Fenchel-Rockefeller duality theorem to obtain a dual formulation.

\begin{lemma}[Convex formulation] \label{lem:convexform} Let $H_{p,q}(\alpha)$ be defined in \eqref{eq:pq3lineproblem}, then for any $r_0,r_1\geq 1$ we have
\begin{align*}
	\sup_{h \in \mathbb{H}^{p,q}(S)} \underbrace{\re h(i\alpha)- \frac{\alpha}{r_1} \norm{h_1}_{L^q}^{r_1} - \frac{1-\alpha}{r_0} \norm{h_0}_{L^p}^{r_0}}_{\coloneqq \mathcal{E}_{p,q,r_0,r_1}^\alpha(h)} = \frac{1}{r_\alpha^\ast} \bigr(H_{p,q}(\alpha)\bigr)^{r_\alpha^\ast},
\end{align*}
where $r_\alpha^{-1} \coloneqq (1-\alpha) r_0^{-1} + \alpha r_1^{-1}$ and $r_\alpha^\ast = \frac{r_\alpha}{r_\alpha-1}$ is the H\"older conjugate of $r_\alpha$, and $\re h(\ii \alpha)$ denotes the real part of $h(\ii \alpha)$. In particular, if a optimizer of Problem~\eqref{eq:pq3lineproblem} exists, then it is unique up to the three-parameter transformation $\tau_{\tau,\beta,\omega}$ defined in~\eqref{eq:symmetry}.
\end{lemma}
\begin{proof} First, from the translation and global phase invariance of $L^p$ norms we find 
\begin{align*}
\sup_{h \in \mathbb{H}^{p,q}(S)} \mathcal{E}_{p,q,r_0,r_1}^\alpha(h) &= \sup_{h} \sup_{\tau >0, |\beta| =1 } \biggr\{ \beta\re h(\ii\alpha -\tau)- \frac{\alpha}{r_1} \norm{h_1}_{L^q}^{r_1} - \frac{1-\alpha}{r_0} \norm{h_0}_{L^p}^{r_0}\biggr\}\\
&= \sup_{h} \biggr\{ \norm{h_\alpha}_{L^\infty}- \frac{\alpha}{r_1} \norm{h_1}_{L^q}^{r_1} - \frac{1-\alpha}{r_0} \norm{h_0}_{L^p}^{r_0}\biggr\}.
\end{align*}
Next, using the transformations $\tau_{0,\alpha \log \omega, \log 1/\omega}$ defined in \eqref{eq:symmetry} we find
\begin{align*}
\sup_{h \in \mathbb{H}^{p,q}(S)} \mathcal{E}_{p,q,r_0,r_1}^\alpha(h) &= \sup_{h \in \mathbb{H}^{p,q}} \sup_{\omega >0} \mathcal{E}_{p,q,r_0,r_1}^\alpha (\tau_{0,\alpha \log \omega,\log 1/\omega} h) \\
&= \sup_{h \in \mathbb{H}^{p,q}(S)} \biggr\{ \norm{h_\alpha}_{L^\infty} - \inf_{\omega >0 } \biggr\{ \underbrace{ \omega^{r_1(\alpha-1)} \frac{\alpha}{r_1} \norm{h_1}_{L^q}^{r_1} + \omega^{s\alpha} \frac{1-\alpha}{r_0} \norm{h_0}_{L^p}^{r_0}}_{\coloneqq f_h(\omega)} \biggr\}  \biggr\}.
\end{align*}
Since $0<\alpha<1$ and $r_0,r_1 >0$, we can minimize the function $f_h(\omega)$ with respect to $\omega>0$ to obtain
\begin{align*}
	\sup_{h \in \mathbb{H}^{p,q}(S)} \mathcal{E}_{p,q,r_0,r_1}^\alpha(h) &=\sup_{h} \biggr\{ \norm{h_\alpha}_{L^\infty} - \frac{1}{r_\alpha} \bigr(\norm{h_0}_{L^p}^{1-\alpha} \norm{h_1}_{L^q}^\alpha\bigr)^{r_\alpha} \biggr\} .
\end{align*}
Finally, by using the transformation $\tau_{0,\beta,0}$ we find
\begin{align*}
	\sup_{h} \mathcal{E}_{p,q,r_0,r_1}^\alpha(h) &= \sup_{h} \sup_{\beta >0} \biggr\{ \underbrace{\beta \norm{h_\alpha}_{L^\infty} - \beta^{r_\alpha}\frac{1}{r_\alpha} \bigr(\norm{h_0}_{L^p}^{1-\alpha} \norm{h_1}_{L^q}^\alpha\bigr)^{r_\alpha}}_{\coloneqq g_h(\beta)} \biggr\},
\end{align*}
Thus since $r_\alpha >1$, we can maximize $g_h(\beta)$ on the interval $\beta>0$ to complete the proof. (The uniqueness statement follows from the strict convexity of $-\mathcal{E}_{p,q,r_0,r_1}$.)
\end{proof}

We now consider $-\mathcal{E}_{p,q,r_0,r_1}^\alpha$ as a convex functional in $L^p(\R)$, and split it as
\begin{align}
	-\mathcal{E}_{p,q,r_0,r_1}^\alpha(f) = \frac{1-\alpha}{r_0} \norm{f}_{L^p}^{r_0} + G(f), \label{eq:splitting}
\end{align}
where $G$ is the functional defined as
\begin{align}
	G(f) = \begin{dcases} \frac{\alpha}{r_1} \norm{h_1}_{L^q}^{r_1} - \re h(\ii\alpha), \quad &\mbox{if $f = h_0$ for some $h \in \mathbb{H}^{p,q}(S)$,} \\
	+\infty, \quad &\mbox{otherwise.} \end{dcases} \label{eq:Gdef}
\end{align}
The Fenchel conjugate of $ \frac{1-\alpha}{r_0} \norm{\cdot}_{L^p}^{r_0}$ in $L^p(\R)$ is well-known and given by
\begin{align*}
	\biggr( \frac{1-\alpha}{r_0} \norm{\cdot}_{L^p}^{r_0}\biggr)^\ast(g) = \frac{(1-\alpha)^{1-r_0^\ast}}{r_0^\ast} \norm{g}_{L^{p^\ast}}^{r_0^\ast}
\end{align*}
where $r_0^\ast = r_0/(r_0-1)$ denotes the H\"older conjugate of $r_0$. The conjugate of $G$ is computed in the next lemma.

\begin{lemma}[Convex conjugate of $G$] \label{lem:dualG} Let $1\leq p < \infty$ and $G:L^p(\R) \rightarrow \R \cup \{+\infty\}$ be the functional defined in \eqref{eq:Gdef}, then its Fenchel conjugate is given by
\begin{align*}
	G^\ast(g) = \begin{dcases} \frac{\alpha^{1-r_1^\ast} }{r_1^\ast} \norm{m_{-1}}_{L^{q^\ast}}^{r_1^\ast} , \quad &\mbox{if $g = m_0$ for some $m$ with $m+\mathscr{P}_\alpha\in \mathbb{H}^{p^\ast,q^\ast}(S^\ast)$,}\\
	+\infty, \quad &\mbox{otherwise,} \end{dcases}
\end{align*}
where $\mathscr{P}_\alpha$ is the Poisson kernel on the strip
\begin{align}
	\mathscr{P}_\alpha(z) \coloneqq \frac{1}{2} \frac{\sin(\pi \alpha)}{\cosh(\pi z) - \cos(\pi \alpha)} ,\quad z\in S.\label{eq:pdef}
\end{align}

\end{lemma}

\begin{proof} Let us assume that the following identity holds and postpone its proof for later.
\begin{align}
	\overline{h(i\alpha)} = \inner{(\mathscr{P}_\alpha)_0, h_0} - \inner{(\mathscr{P}_\alpha)_{-1}, h_1}, \quad \mbox{for any $h \in \mathbb{H}^{p,q}(S)$.} \label{eq:claim}
\end{align}
Then from the definition of Fenchel conjugate we have, for $g\in L^{p^\ast}(\R)$,
\begin{align*}
	G^\ast(g) &= \sup_{h \in \mathbb{H}^{p,q}(S)} \biggr\{ \re \inner{g, h_0} + \re h(i\alpha) - \frac{\alpha}{r_1} \norm{h_1}_{L^q}^{r_1} \biggr\} \\
	&= \sup_{h \in \mathbb{H}^{p,q}(S)} \biggr\{ \re \inner{g + (\mathscr{P}_\alpha)_0, h_0} - \re \inner{(\mathscr{P}_\alpha)_{-1},h_1} - \frac{\alpha}{r_1} \norm{h_1}_{L^q}^{r_1} \biggr\}  \\
	&= \sup_{h \in \mathbb{H}^{p,q}(S)} \sup_{\kappa >0} \biggr\{ \underbrace{\kappa \, \re \biggr(\inner{g+(\mathscr{P}_\alpha)_0,h_0} - \inner{(\mathscr{P}_\alpha)_{-1},h_1}\biggr) -  \kappa^{r_1} \frac{\alpha}{r_1} \norm{h_1}_{L^q}^{r_1}}_{\coloneqq f_h(\kappa)} \biggr\}.
\end{align*}
Since $r_1 >1$, we can maximize the function $f_h(\kappa)$ defined above to obtain
\begin{align*}
	G^\ast(g) =  \frac{\alpha^{1-r_1^\ast}}{r_1^\ast} \biggr(\sup_{h \in \mathbb{H}^{p,q}(S)} \biggr\{\re \frac{\inner{g+(\mathscr{P}_\alpha)_0,h_0} - \inner{(\mathscr{P}_\alpha)_{-1},h_1}}{\norm{h_1}_{L^q}}\biggr\}\biggr)^{r_1^\ast}.
\end{align*}
Since $(\mathscr{P}_\alpha)_{-1} \in L^{q^\ast}(\R)$, the supremum on the right hand side is bounded if and only if 
\begin{align*}
	|\inner{g+(\mathscr{P}_\alpha)_0, h_0}| \lesssim \norm{h_1}_{L^q} \mbox{\quad for any $h \in \mathbb{H}^{p,q}(S)$.}
\end{align*}
Lemma~\ref{lem:boundarycriteria} then implies that $g \in \dom G^\ast$ if and only if $g = m_0$ for some $m+\mathscr{P}_\alpha \in  \mathbb{H}^{p^\ast, q^\ast}(S^\ast)$. Moreover, in this case we have (see eq.~\eqref{eq:closedcurveintegral})
\begin{align*}
G^\ast(g) = \frac{\alpha^{1-r_1^\ast}}{r_1^\ast}  \biggr(\sup_{h \in \mathbb{H}^{p,q}(S)} \biggr\{\re \frac{\inner{(m+\mathscr{P}_\alpha)_{-1},h_1} - \inner{(\mathscr{P}_\alpha)_{-1},h_1}}{\norm{h_1}_{L^q}}\biggr\}\biggr)^{r_1^\ast} = \frac{\alpha^{1-r_1^\ast}}{r_1^\ast} \norm{m_{-1}}_{L^{q^\ast}}^{r_1^\ast},
\end{align*}
which proves the lemma. 

Let us now prove eq.~\eqref{eq:claim}. To this end, note that the function $F(z) = \mathscr{P}_\alpha(z) h^\ast(z)$ is meromorphic on $S^\ast$ with a simple pole at $-i\alpha$. Therefore, from Cauchy's residue theorem we have 
\begin{align*}
	\overline{h(i\alpha)} = -(2\pi i) \lim_{z \ra -i\alpha} F(z)(z+i\alpha)  = \int_{-R}^R \bigr( F_{-\epsilon}(x)-F_{\epsilon-1}(x)\bigr) \mathrm{d}x + \int_{\epsilon-1}^{-\epsilon} \bigr( F(-R + iy) - F(R + iy) \bigr)\mathrm{d}y 
\end{align*}
for any $R>0$ and $\epsilon>0$. But since $F(x+iy) \ra 0$ as $|x| \ra 0$ uniformly for $y$ in compact subsets of $(0,1)$, by taking the limit $R\ra \infty$ we find that
\begin{align*}
	\inner{(\mathscr{P}_\alpha)_{-\epsilon}, h_{\epsilon}} - \inner{(\mathscr{P}_\alpha)_{\epsilon-1}, h_{1-\epsilon}} = \int_\R F_{-\epsilon}(x) \mathrm{d} x - \int_\R F_{\epsilon-1}(x) \mathrm{d}x =  \overline{h(i\alpha)}.
\end{align*}
We can now pass to the limit $\epsilon \ra 0^+$ to prove~\eqref{eq:claim} because $(\mathscr{P}_\alpha)_{-\epsilon}, (\mathscr{P}_\alpha)_{\epsilon-1} \ra (\mathscr{P}_\alpha)_0, (\mathscr{P}_\alpha)_{-1}$ strongly in $L^1(\R)\cap L^\infty(\R)$ and $h_{\epsilon},h_{1-\epsilon} \ra h_0, h_1$ weakly in $L^p(\R)$, respectively $L^q(\R)$ (or weakly-$\star$ in the case $p$ or $q=\infty$), as $\epsilon \ra 0^+$.
\end{proof}

\begin{remark*} The Poisson kernel $\mathscr{P}_\alpha$ is not specially important here; it could be replaced for any other holomorphic function on $S^\ast\setminus \{-\ii \alpha\}$ with a simple pole with residue $i/2\pi$ at $-\ii\alpha$ and decaying fast enough away from this pole.
\end{remark*}

We are now in position to present the main result of this subsection, namely, a dual formulation of Problem~\eqref{eq:pq3lineproblem} in a suitable space of meromorphic functions.

\begin{theorem}[Dual problem - meromorphic version] \label{thm:dual} Let $1\leq p < \infty, 1\leq q \leq\infty$, $H_{p,q}(\alpha)$ be defined via eq.~\eqref{eq:pq3lineproblem}, and $\mathscr{P}_\alpha$ be the Poisson kernel defined in~\eqref{eq:pdef}. Then we have
\begin{align}
	H_{p,q}(\alpha) =  \min \biggr\{\biggr(\frac{\norm{m_0}_{L^{p^\ast}}}{1-\alpha}\biggr)^{1-\alpha} \biggr(\frac{\norm{m_{-1}}_{L^{q^\ast}}}{\alpha}\biggr)^\alpha : m \in \mathscr{P}_\alpha + \mathbb{H}^{p^\ast,q^\ast}(S^\ast)\biggr\}, \label{eq:dualproblem}
\end{align}
where $\frac{1}{p^\ast} + \frac{1}{p} = 1 = \frac{1}{q} + \frac{1}{q^\ast}$. Moreover, the minimizer in \eqref{eq:dualproblem} exists and is unique up to the transformation $\tau_{\kappa} m = m(x) \kappa^{iz-\alpha} (\kappa >0)$.
\end{theorem}

\begin{proof} First, note that $f \mapsto \norm{f}_{L^p}^p$ is continuous on $L^p$ and $\dom G  \cap L^p(\R) \neq \emptyset$. So the sum defined in \eqref{eq:splitting} satisfies the assumptions on the Fenchel-Rockafellar theorem \cite{Roc72}. Consequently, by the Fenchel-Rockafellar duality theorem, Lemma~\ref{lem:convexform} and Lemma~\ref{lem:dualG} we find (see~\eqref{eq:splitting})
\begin{align*}
	\frac{1}{r_\alpha^\ast} \bigr(H_{p,q}(\alpha)\bigr)^{r_\alpha^\ast} &= \min \biggr\{ \frac{(1-\alpha)^{1-r_0^\ast}}{r_0^\ast} \norm{g}_{L^{p^\ast}}^{r_0^\ast} + G^\ast(-g): g\in L^{p^\ast}(\R) \biggr\}\\
 &=\min \biggr\{ \frac{(1-\alpha)^{1-r_0^\ast}}{r_0^\ast} \norm{m_0}_{L^{p^\ast}}^{r_0^\ast}+ \frac{\alpha^{1-r_1^\ast}}{r_1^\ast} \norm{m_{-1}}_{L^{q^\ast}}^{r_1^\ast} : m \in \mathscr{P}_\alpha + \mathbb{H}^{p^\ast,q^\ast}(S^\ast)\biggr\},
\end{align*}
where the minimizer exists (and is unique by the strict convexity of $\norm{\cdot}^{s}$ with $s>1$). Next, note that the function
\begin{align*} \mathscr{P}_\alpha(z) \kappa^{iz -\alpha} - \mathscr{P}_\alpha(z) \quad \mbox{ belongs to $\mathbb{H}^{p,q}(S^\ast)$ for any $1\leq p,q\leq \infty$.}
\end{align*}
Indeed, since $\kappa^{i(-i\alpha)-\alpha} = 1$, the poles of $\mathscr{P}_\alpha$ and $\kappa^{\ii z-\alpha} \mathscr{P}_\alpha$ cancel out while the decay as $|x| \rightarrow \infty$ remains. In particular,
\begin{align*}
	\tau_\kappa m(z) = m(z) \kappa^{iz - \alpha} \in \mathscr{P}_\alpha + \mathbb{H}^{p^\ast,q^\ast}(S^\ast), \quad\mbox{for any $m \in \mathscr{P}_\alpha + \mathbb{H}^{p^\ast,q^\ast}(S^\ast)$ and $\kappa>0$.}
\end{align*}
Therefore,
\begin{align*}
	\frac{1}{r_\alpha^\ast} \bigr(H_{p,q}(\alpha)\bigr)^{r_\alpha^\ast} &= \min_{m} \inf_{\kappa >0}\biggr\{ \frac{(1-\alpha)^{1-r_0^\ast}}{r_0^\ast} \norm{(\tau_\kappa m)_0}_{L^{p^\ast}}^{r_0^\ast} + \frac{\alpha^{1-r_1^\ast}}{r_1^\ast} \norm{(\tau_\kappa m)_{-1}}_{L^{q^\ast}}^{r_1^\ast} \biggr\} \\
	&=  \min_{m} \inf_{\kappa>0} \biggr\{ \kappa^{-r_0^\ast \alpha} \frac{1-\alpha}{r_0^\ast}\biggr(\frac{\norm{m_0}_{L^{p^\ast}}}{1-\alpha}\biggr)^{r_0^\ast} + \kappa^{r_1^\ast(1-\alpha)} \frac{\alpha}{r_1^\ast} \biggr(\frac{\norm{m_{-1}}_{L^{q^\ast}}}{\alpha}\biggr)^{r_1^\ast} \biggr\} \\
	&=   \frac{1}{r_\alpha^\ast} \biggr(\min_{m} \biggr\{ \biggr( \frac{\norm{m_0}_{L^{p^\ast}}}{1-\alpha}\biggr)^{1-\alpha} \biggr(\frac{\norm{m_{-1}}_{L^{q^\ast}}}{\alpha}\biggr)^\alpha\biggr\}\biggr)^{r_\alpha^\ast},
\end{align*}
which completes the proof.
\end{proof}

Before proceeding, let us make a few comments regarding the restriction $p<\infty$ in Theorem~\ref{thm:dual}. First, we note that $H_{p,q}(\alpha)=H_{q,p}(1-\alpha)$ by flipping the strip. Hence, it is enough to consider the case $p\leq q$ so that the restriction $p<\infty$ is only relevant for the case $p=q=\infty$. For this case, one could in fact introduce the space of finitely additive measures that are absolutely continuous with respect to Lebesgue measure, i.e., the dual of $L^\infty(\R)$, and follow the same approach to obtain the corresponding results. However, this requires additional technical complications that we seek to avoid here. 

Fortunately, as we shall show next, the dual problem~\eqref{eq:dualproblem} is equivalent to the primal problem associated to $H_{p^\ast,q^\ast}(\alpha)$, where $p^\ast,q^\ast$ are the H\"older conjugate exponents of $p$ and $q$. Therefore, we obtain the duality relation in Theorem~\ref{thm:dual} for $p= q=\infty$ by considering the case $p=q=1$. In this way, however, we do not immediately obtain the existence of optimizers for problem~\eqref{eq:pq3lineproblem} in the case $p=q=1$ (as their existence was a consequence of the Fenchel duality theorem employed in the proof of Theorem~\ref{thm:dual}). Since we anyways construct the optimizers explicitly in Section~\ref{sec:eulerlagrange}, this observation will play no role in our subsequent analysis.

\subsection{Duality relation}

We now turn to the proof of Theorem~\ref{thm:duality}. The key idea here is to use a factorization of the space $\mathscr{P}_\alpha+ \mathbb{H}^{p^\ast,q^\ast}(S^\ast)$ to show that the dual problem~\eqref{eq:dualproblem} is equivalent to the three lines problem~\eqref{eq:pq3lineproblem} in $\mathbb{H}^{p^\ast,q^\ast}(S)$. To state this factorization precisely, let us define the Blaschke factor on the strip $S$ with zero at $\ii \alpha$ as
\begin{align}
    B_\alpha(z) \coloneqq \ee^{-\ii \pi \alpha} \frac{\ee^{\ii \pi \alpha}-\ee^{\pi z}}{\ee^{\pi z} - \ee^{-\ii \pi \alpha}}. \label{eq:Blaschkedef}
\end{align}
Then the following result holds.
    \begin{lemma}[Blaschke factorization]\label{lem:blaschkefactorization} Let $\alpha \in (0,1)$ and $1\leq p,q \leq \infty$, and let $\mathscr{P}_\alpha$ be the Poisson kernel defined in~\eqref{eq:pdef}. Then $m -\mathscr{P}_\alpha \in \mathbb{H}^{p,q}(S^\ast)$ if and only if the function $g \in \mathbb{H}(S)$ defined as
    \begin{align}
        g(z) = m^\ast(z) B_\alpha(z) \label{eq:tildemdef}
    \end{align}
    belongs to $\mathbb{H}^{p,q}(S)$ and satisfies 
    \begin{align}
        g(\ii \alpha) = \frac{1}{4 \sin(\pi \alpha)}. \label{eq:gvalue}
    \end{align}
    Moreover, in this case we have
    \begin{align}
        \norm{g_0}_{L^{p}(\R)} = \norm{m_0}_{L^{p}(\R)} \quad \mbox{and}\quad \norm{g_1}_{L^{q}(\R)} = \norm{m_{-1}}_{L^{q}(\R)}. \label{eq:normids}
    \end{align}
\end{lemma}

\begin{proof} Note that the function $B_\alpha$ is bounded inside $S$ and has a single simple zero at $\ii \alpha$. In particular, the function $z \mapsto \mathscr{P}_\alpha^\ast(z) B_\alpha(z) = \mathscr{P}_\alpha(z) B_\alpha(z)$ belongs to $\mathbb{H}^{p,q}(S)$ for any $1\leq p,q\leq \infty$. Similarly, $m B_\alpha \in \mathbb{H}^{p,q}(S)$ for any $m \in \mathbb{H}^{p,q}(S)$. Consequently the function $g$ defined in~\eqref{eq:tildemdef} satisfies 
\begin{align*}
    g(z) = (m-\mathscr{P}_\alpha)^\ast(z) B_\alpha(z) + \mathscr{P}_\alpha(z) B_\alpha(z) \in \mathbb{H}^{p^\ast,q^\ast}(S). 
\end{align*}
The identities in \eqref{eq:normids} and \eqref{eq:gvalue} follow respectively from the facts that $|B_\alpha(z)| = 1$ for any $z \in \partial S$ and 
\begin{align*}
    \lim_{z \rightarrow \ii \alpha} \mathscr{P}_\alpha(z) B_\alpha(z) = \frac{1}{4 \sin(\pi \alpha)}.
\end{align*}
The converse follows by noticing that $\frac{1}{B_\alpha(z)}$ is uniformly bounded away from $z=\ii \alpha$, and therefore, $g/B_\alpha - \mathscr{P}_\alpha \in \mathbb{H}^{p,q}(S)$ provided that $g \in \mathbb{H}^{p,q}(S)$ and \eqref{eq:gvalue} holds.
\end{proof}

We can now complete the proof of Theorem~\ref{thm:duality}.

\begin{proof}[Proof of Theorem~\ref{thm:duality}]
First, note that the transformation $m \mapsto g$ defined via \eqref{eq:tildemdef} gives an one-to-one mapping between the domain of the dual problem~\eqref{eq:dualproblem} and the set
\begin{align*}
    \mathcal{D} \coloneqq \biggr\{ g \in \mathbb{H}^{p^\ast,q^\ast}(S) : g(\ii \alpha) = \frac{1}{4 \sin(\pi \alpha)} \biggr\}.
\end{align*}
We can thus restate the dual problem in~\eqref{eq:dualproblem} as an equivalent problem over the set $\mathcal{D}$. Hence, from the duality in Theorem~\ref{thm:dual} we have
\begin{align}
    H_{p,q}(\alpha) &= \frac{1}{\alpha^\alpha (1-\alpha)^{1-\alpha}} \min \biggr\{ \frac{\norm{g_0}_{L^{p^\ast}}^{1-\alpha} \norm{g_1}_{L^{q^\ast}}^\alpha}{4 \sin(\pi \alpha) |g(\ii \alpha)|} : g \in \mathcal{D} \biggr\}\nonumber \\
    &= \frac{1}{4 \sin(\pi \alpha) \alpha^\alpha (1-\alpha)^{1-\alpha}} \biggr(\max \biggr\{ \frac{|g(\ii \alpha)|}{\norm{g_0}_{L^{p^\ast}}^{1-\alpha} \norm{g_1}_{L^{q^\ast}}^\alpha} : g\in \mathcal{D}\biggr\}\biggr)^{-1} \label{eq:middlemax}
\end{align}
Now note that, due to the invariance of $|g(\ii \alpha)|/(\norm{g_0}_{L^{p^\ast}}^{1-\alpha} \norm{g_1}_{L^{q^\ast}}^\alpha)$ under scalar multiplication $g\mapsto \beta g$ with $\beta \in \C\setminus \{0\}$, the maximization problem in \eqref{eq:middlemax} can be relaxed to the whole space $\mathbb{H}^{p^\ast,q^\ast}(S)$. Moreover, by the translation invariance of the $L^p$ norms, we can replace $|g(\ii \alpha)|$ in \eqref{eq:middlemax} by $\norm{g}_{L^\infty}$. We therefore conclude that
\begin{align*}
    \frac{1}{4 \sin(\pi \alpha) \alpha^\alpha (1-\alpha)^{1-\alpha}} \biggr(\max \biggr\{ \frac{|g(\ii \alpha)|}{\norm{g_0}_{L^{p^\ast}}^{1-\alpha} \norm{g_1}_{L^{q^\ast}}^\alpha} : g\in \mathcal{D}\biggr\}\biggr)^{-1}  = \frac{1}{4 \sin(\pi \alpha) \alpha^\alpha (1-\alpha)^{1-\alpha}} \frac{1}{H_{p^\ast,q^\ast}(\alpha)},
\end{align*}
which together with~\eqref{eq:middlemax} completes the proof of \eqref{eq:dualityrelation}. The existence and uniqueness of maximizers of $H_{p^\ast,q^\ast}(S)$ follows from the corresponding results (Theorem~\ref{thm:dual}) for the dual problem~\eqref{eq:dualproblem}.
\end{proof}
\section{The Euler-Lagrange equations}
\label{sec:eulerlagrange}

In this section we effectively solve the Euler-Lagrange equations associated with the primal and dual problems in \eqref{eq:pq3lineproblem} and \eqref{eq:dualproblem}.

Let us start by writing down the E.L. equation for the dual problem~\eqref{eq:dualproblem}. In this case, if we denote by
\begin{align}
    \mathcal{E}^\ast_{p^\ast,q^\ast,\alpha}(m) \coloneqq  \biggr(\frac{\norm{m_0}_{L^{p^\ast}}}{1-\alpha}\biggr)^{1-\alpha} \biggr(\frac{\norm{m_{-1}}_{L^{q^\ast}}}{\alpha}\biggr)^\alpha, \quad \dom \mathcal{E}^\ast_{p^\ast,q^\ast,\alpha} = \mathscr{P}_\alpha + \mathbb{H}^{p^\ast,q^\ast}(S^\ast)
\end{align}
then straightforward computation shows that the Gateaux derivative of $\mathcal{E}_{p^\ast,q^\ast}^\ast$ at $m \in \dom \mathcal{E}_{p^\ast,q^\ast}^\ast$ is given by 
\begin{align*}
    &\mathrm{d}_m \mathcal{E}^\ast_{p^\ast,q^\ast,\alpha}(\phi) = \frac{c^\alpha}{\norm{m_0}_{L^{p^\ast}}^{p^\ast-1}} \re \inner{ m_0 |m_0|^{p^\ast-2}, \phi_0} + \frac{c^{\alpha-1}}{\norm{m_{-1}}_{L^{q^\ast}}^{q^\ast-1}} \re \inner{m_{-1} |m_{-1}|^{q^\ast-2}, \phi_{-1}}
\end{align*}
for any $\phi \in \mathbb{H}^{p^\ast,q^\ast}(S^\ast)$, where
\begin{align}
   c = \frac{(1-\alpha) \norm{m_{-1}}_{L^{q^\ast}}}{\alpha \norm{m_0}_{L^{p^\ast}}}. \label{eq:cdef}
\end{align}
The above Euler-Lagrange equations are rigorous except for the cases where either $p \in \{1,\infty\}$ or $q \in \{1,\infty\}$. Indeed, in the case $p=\infty$ (respectively $q=\infty$), the Gateaux derivative of $\norm{\cdot}_{L^p}$ is only well-defined if $m_0(x)$ (respectively $m_{-1}(x)$) is almost everywhere non-vanishing\footnote{This difficulty can be circumvented by using the product Lemma~\ref{lem:product} and the well-known fact (see, e.g., \cite[Corollary 4.2]{Gar06}) that the trace of a non-zero function in the classical Hardy space $\mathbb{H}^1(S)$ can not vanish in a subset of positive measure in $\R$. We shall not, however, use this fact here.}. The difficulty in the case $p=1$ (or $q=1$) is even more extreme because $p^\ast = \infty$ and we can give no meaningful sense to the function $m_0|m_0|^{p^\ast-2}$. To overcome these issues, we shall instead consider the primal-dual pair $(h,m) \in \mathbb{H}^{p,q}(S) \times (\mathscr{P}_\alpha + \mathbb{H}^{p^\ast,q^\ast}(S^\ast))$ and derive a joint Euler-Lagrange equation for them. 

\begin{lemma}[Euler-Lagrange equations] \label{lem:eulerlagrange} Let $\alpha\in (0,1)$ and $1\leq p,q \leq \infty$ and suppose that $h \in \mathbb{H}^{p,q}(S)$ and $m \in \mathscr{P}_\alpha + \mathbb{H}^{p^\ast,q^\ast}(S^\ast)$ satisfy 
\begin{align}
    h_0 = m_0 |m_0|^{p^\ast-2} \quad \mbox{or}\quad m_0 = h_0 |h_0|^{p-2}  \label{eq:ELeq0}
    \intertext{and}
    h_1 = - c m_{-1}|m_{-1}|^{q^\ast-2} \quad\mbox{or} \quad m_{-1} = -c^{1-q} h_1 |h_1|^{q-2}, \label{eq:ELeq1}
\end{align}
for some constant $c>0$. Then $h$ is an optimizer of problem~\eqref{eq:pq3lineproblem} and $m$ is an optimizer of the dual problem~\eqref{eq:dualproblem}. 
\end{lemma}

Let us now make some comments on the Euler-Lagrange equations~\eqref{eq:ELeq0} and \eqref{eq:ELeq1}. First, we note that in the special case where either $p=2$ or $q=2$, equations~\eqref{eq:ELeq0} and \eqref{eq:ELeq1} imply that $m$ and $h$ are in fact meromorphic/holomorphic extensions of each other. In \cite{CCR24}, this observation (together with a different Blaschke factorization as the one stated in Lemma~\ref{lem:blaschkefactorization}) allowed to reduce the problem to a Euler-Lagrange equation for a single holomorphic function, for which a solution could be constructed via the Fourier transform of a principal value distribution \cite[Lemma 4.6]{CCR24}. Extending this approach to the general $p,q$ case is a non-trivial task mostly because it requires the construction of two phase functions (one for $h$ and one for $m$) satisfying (linear versions) of equations~\eqref{eq:ELeq0} and \eqref{eq:ELeq1}.

For the proof of Lemma~\ref{lem:eulerlagrange}, we shall need an additional lemma. This lemma turns out to considerably simplify the approach used in \cite[Section 4]{CCR24}; it gives an explicit formula for the product of primal and dual optimizers by combining equations~\eqref{eq:ELeq0} and \eqref{eq:ELeq1} with the Schwarz reflection principle (Lemma~\ref{lem:schwarzreflection}) and Liouville's theorem.

\begin{lemma}[Dual-primal product formula]\label{lem:dualproduct} Let $p,q\in [1,\infty]$ and $h \in \mathbb{H}^{p,q}(S)$ and $m \in \mathscr{P}_\alpha +\mathbb{H}^{p^\ast,q^\ast}(S^\ast)$ be such that \eqref{eq:ELeq0} and \eqref{eq:ELeq1} holds. Then we have
\begin{align}
    m^\ast(z) h(z) = \frac{\kappa}{1-\alpha} \mathscr{P}_\alpha(z), \quad \mbox{for any $z\in S$,} \label{eq:periodicidentity}
\end{align}
where 
\begin{align}
    \kappa = \begin{dcases} \norm{m_0}_{L^{p^\ast}}^{p^\ast}, \quad &\mbox{if $p^\ast \neq \infty$}\\
    \norm{h_0}_{L^p}^p, \quad &\mbox{if $p\neq \infty$.} \end{dcases} \label{eq:alphadef}
    \end{align}
\end{lemma}

\begin{proof} Let $F(z) \coloneqq h(z) m^\ast(z)$. Since $h$ is holomorphic on $S$ and $m$ is meromorphic on $S^\ast$ with a simple pole at $-\ii \alpha$, the function $F$ is meromorphic on $S$ and has at most a single simple pole at $\ii \alpha$. Moreover, $F$  belongs to $\mathbb{H}^{1,\infty}(S_{0,a})$ for any $0<a<\alpha$ and satisfies
    \begin{align}
    F(x) = |m_0(x)|^{p^\ast} \quad \mbox{or} \quad F(x) = |h_0(x)|^{p} \quad \mbox{(or both in case $1<p<\infty$)}  \label{eq:reflectionat0} \end{align}
    for almost every $x\in \R$. Hence, by the Schwarz reflection principle, $F$ can be holomorphically extended to a function in $S_{-a,0}$ satisfying
    \begin{align}
        F(z) = F^\ast(z), \quad \mbox{for $z \in S_{-a,a}$.} \label{eq:reflectedF}
    \end{align}
    Thus by analytic continuation, $F$ can be meromorphically extended further to a function in $S_{-1,1}$ satisfying \eqref{eq:reflectedF} for $z\in S_{-1,1}$. In particular, this extension (denoted again by $F$) has two simple poles around the points $z = \ii \alpha$ and $z = -\ii \alpha$ and satisfies
    \begin{align}
        \beta \coloneqq \lim_{z \ra \ii \alpha} F(z) (z-\ii \alpha) = \lim_{z\ra -\ii \alpha} \overline{F(z)(z+\ii \alpha)}. \label{eq:betapole}
    \end{align}
    
    Next, note that $F$ also satisfies $F \in \mathbb{H}^{\infty,1}(S_{a,1})$ for any $\alpha < a<1$ and
    \begin{align*}
        F_1(x) = -c |m_{-1}|^{q^\ast} \quad\mbox{or}\quad F_1(x) = -c^{q-1} |h_1(x)|^q \quad \mbox{(or both)}
    \end{align*}
for almost every $x\in \R$ and some $c>0$. Thus by applying the Schwarz relfection principle to the function $G(z) \coloneqq F(z+\ii)$, we can further meromorphically extend $F$ to $S_{-1,3}$ by setting
\begin{align*}
    F(z) = \overline{F(\overline{z}+\ii)}, \quad \mbox{$z\in S_{1,3}$}.
\end{align*}
In summary, the function $F$ is now a meromorphic function in $S_{-1,3}$ with simple poles at the points $\{-\ii \alpha, \ii \alpha, \ii (2-\alpha)m, \ii(2+\alpha)\}$, and satisfying $F_{-1}(x) = F_1(x) = F_3(x) \in \R$ and $F_0(x) = F_2(x) \in \R$ for almost every $x\in \R$. Therefore, repeating the preceding arguments, we can meromorphically extend $F$ to the whole complex plane by successively reflecting $F$ over the lines $\R + \ii k$, $k \in \Z$. This extension then satisfies $F(z+2\ii) = F(z)$ ($2\ii$-periodic) for any $z\in \C$, and its poles are all simple and located at $\{\pm \ii \alpha + 2\ii k : k \in \Z\}$. 

We now claim that this is enough to conclude that
\begin{align}
    F(z) = \frac{\kappa}{1-\alpha} \mathscr{P}_\alpha(z) \quad \mbox{with $\kappa$ given by \eqref{eq:alphadef}.}
\end{align}
Indeed, let $\beta \in \C$ be defined by \eqref{eq:betapole} and consider the function
\begin{align}
    G(z) \coloneqq F(z) + 2 \pi \mathrm{Im}\, \beta \mathscr{P}_\alpha(z) - 2\pi \mathrm{Re}\, \beta \frac{\sinh(\pi z)}{\sin(\pi \alpha)} \mathscr{P}_\alpha(z). 
\end{align}
Then $G$ is meromorphic and $2\ii$-periodic in $\C$, and its poles are at most simple and contained in the set $\{ \pm \ii \alpha + \ii 2k : k \in \Z\}$. However, by construction we have
\begin{align*}
    \lim_{z\ra \ii\alpha} G(z)(z-\ii \alpha) &= \beta + 2 \pi \mathrm{Im}\, \beta \lim_{z\ra \ii \alpha} (z-\ii \alpha) \mathscr{P}_\alpha(z) - 2 \pi \re \beta \lim_{z\ra \ii \alpha} (z-\ii \alpha) \frac{\sinh(\pi z)}{\sin(\pi\alpha)} \mathscr{P}_\alpha(z) \\
    &= \beta + (2 \pi \mathrm{Im}\, \beta - \ii 2 \pi \re \beta) \lim_{z\ra \ii \alpha} (z-\ii\alpha) \mathscr{P}_\alpha(z)  = 0.
    \intertext{and}
    \lim_{z\ra -\ii\alpha} G(z)(z+\ii \alpha) &= \overline{\beta} + 2 \pi \mathrm{Im}\, \beta \lim_{z\ra -\ii \alpha} (z+\ii \alpha) \mathscr{P}_\alpha(z) - 2 \pi \re \beta \lim_{z\ra -\ii \alpha} (z+\ii \alpha) \frac{\sinh(\pi z)}{\sin(\pi \alpha)} \mathscr{P}_\alpha(z)\\
    & \overline{\beta} + (2 \pi \mathrm{Im}\, \beta + \ii 2 \pi \re \beta) \lim_{z\ra \ii \alpha} (z+\ii\alpha) \mathscr{P}_\alpha(z)  = 0,
\end{align*}
which implies that $G$ is in fact entire and $2\ii$-periodic. Since $F(z)$ is bouded away from its poles, we conclude that $G$ is bounded and entire so that by Liouville's theorem, $G = c$ for some constant $c\in \C$. To complete the proof, just note that $c = 0 = \re \beta$ because $G_0$ is integrable along $\R$, and $\mathrm{Im} \beta \neq 0$ because $F \neq 0$. Hence \eqref{eq:periodicidentity} holds and the formula for $\kappa$ follows from~\eqref{eq:reflectionat0}.
\end{proof}

We now turn to the proof of Lemma~\ref{lem:eulerlagrange}.
\begin{proof}[Proof of Lemma~\ref{lem:eulerlagrange}]
Under the assumption of Lemma~\ref{lem:eulerlagrange}, we have an explicit formula for the product $F(z) = h(z) \overline{m(\overline{z})}$ (by Lemma~\ref{lem:dualproduct}). The idea is then to use this formula to prove that the duality gap between $m$ and $h$ is zero. For this, first note that, since $m \in \mathscr{P}_\alpha + \mathbb{H}^{p^\ast,q^\ast}(S^\ast)$, we have
\begin{align}
    \mathrm{Res}(F, \ii \alpha) = \frac{\kappa}{1-\alpha} = h(\ii \alpha).  \label{eq:residue1}
\end{align}
Now note that from the formula $\int_\R \mathscr{P}_{1-\alpha}(x) \mathrm{d} x = \alpha$ (see Lemma~\ref{lem:PoissonFourier}) and \eqref{eq:ELeq1}, we have
\begin{align*}
   -\int_{\R} F_1(x) \mathrm{d} x = \frac{\alpha \kappa}{1-\alpha} = \begin{dcases} c \norm{m_{-1}}_{L^{q^\ast}}^{q^\ast}, \quad &\mbox{if $q^\ast \neq \infty$,} \\
   c^{1-q} \norm{h_1}_{L^q}^q, &\mbox{if $q \neq \infty$.} \end{dcases}
\end{align*}
Either way, we have (even when $q=\infty$ or $q=1$)
\begin{align}
    \norm{h_1}_{L^q} = \biggr(\frac{\alpha \kappa}{1-\alpha}\biggr)^{\frac{1}{q}} c^{\frac{1}{q^\ast}} \quad \mbox{and} \quad \norm{m_{-1}}_{L^{q^\ast}} = \biggr(\frac{\alpha \kappa}{1-\alpha}\biggr)^{\frac{1}{q^\ast}} \frac{1}{c^{\frac{1}{q^\ast}}}. \label{eq:normat1}
\end{align}
Similarly, from \eqref{eq:ELeq0} we find (even when $p=1$ or $p=\infty$),
\begin{align}
    \norm{h_0}_{L^p} = \kappa^{\frac{1}{p}} \quad \mbox{and} \quad \norm{m_0}_{L^{p^\ast}} = \kappa^{\frac{1}{p^\ast}}. \label{eq:normat0}
\end{align}
Thus from \eqref{eq:residue1}, \eqref{eq:normat1}, and~\eqref{eq:normat0}, we conclude that
\begin{align*}
    \frac{\norm{h_\alpha}_{L^\infty}}{\norm{h_0}_{L^p}^{1-\alpha} \norm{h_1}_{L^q}^\alpha} \geq \frac{h(\ii \alpha)}{\norm{h_0}_{L^p}^{1-\alpha} \norm{h_1}_{L^q}^\alpha} = \frac{\kappa^{\frac{1-\alpha}{p^\ast} + \frac{\alpha}{q^\ast}}}{1-\alpha} \biggr(\frac{1-\alpha}{\alpha}\biggr)^{\frac{\alpha}{q}} \frac{1}{c^{\frac{\alpha}{q^\ast}}} = \frac{\norm{m_0}_{L^{p^\ast}}^{1-\alpha} \norm{m_{-1}}_{L^{q^\ast}}^\alpha}{(1-\alpha)^{1-\alpha} \alpha^\alpha}. 
\end{align*}
The result now follows from Theorem~\ref{thm:dual}.
\end{proof}

As previously mentioned, Lemma~\ref{lem:dualproduct} allows us to use the Poisson representation formula in the strip (Lemma~\ref{lem:PoissonRepresentation}) to construct the phase functions corresponding to the solutions of the Euler-Lagrange equations~\eqref{eq:ELeq0} and \eqref{eq:ELeq1}. For this, we shall use the following lemma.

\begin{lemma}[Construction of optimizers] \label{lem:phasefunction} Let $0 < \alpha < 1$ and $1\leq p,q\leq \infty$, then there exists a unique holomorphic function $\phi_{\alpha,p,q}: \mathcal{S} \rightarrow \C$ satisfying $\mathrm{Im}\, \phi_{\alpha,p,q}(\ii \alpha) = 0$ and 
\begin{align}
    \re \, \phi_{\alpha,p,q}(x+\ii y) = \frac{1}{p} \mathscr{P}_{y} \ast \log\left(\frac{\mathscr{P}_\alpha}{1-\alpha}\right) + \frac{1}{q} \mathscr{P}_{1-y} \ast \log\left(\frac{\mathscr{P}_{1-\alpha}}{\alpha}\right). \label{eq:optimizerformula}
\end{align}
Moreover, we have
\begin{align}
    \sup_{x + \ii y \in \mathcal{S}} \left|\re \phi_{\alpha,p,q}(x+\ii y) + \pi \biggr(\frac{1-y}{p}  + \frac{y}{q} \biggr) |x|\right| < \infty. \label{eq:phasedecay}
\end{align}
\end{lemma}

\begin{proof} Note that the function $\re \phi_{\alpha,p,q}$ defined via \eqref{eq:optimizerformula} is harmonic. As $S$ is simply connected, the uniqueness of $\phi_{\alpha,p,q}$ follows from the existence of an unique (up to a constant) harmonic conjugate on $S$ and the constraint $\mathrm{Im}\, \phi_{\alpha,p,q}(\ii \alpha) =0$. To prove~\eqref{eq:phasedecay} we note that 
\begin{align*}
    1\lesssim_\alpha \mathscr{P}_\alpha(x) \ee^{\pi |x|} \lesssim_\alpha 1
\end{align*}
with an implicit constant independent of $x\in \R$, and
\begin{align*}
    \int_{\R} \mathscr{P}_y(w) (1+|w|) \mathrm{d} w \lesssim 1
\end{align*}
with an implicit constant independent of $0<y<1$. Thus from Lemma~\ref{lem:PoissonFourier}, we have
\begin{align*}
    \re \phi_{\alpha,p,q}(x+\ii y) + \pi \biggr(\frac{1-y}{p} + \frac{y}{q}\biggr)|x| =& \int_\R \frac{\mathscr{P}_y(w)}{p} \biggr(\log\left(\frac{\mathscr{P}_\alpha(x-w)}{1-\alpha}\right)+ \log(\ee^{\pi |x|})\biggr) \mathrm{d} w \\ 
    &+ \int_\R \frac{\mathscr{P}_{1-y}(w)}{q} \biggr(\log\left(\frac{\mathscr{P}_{1-\alpha}(x-w)}{\alpha}\right)+ \log(\ee^{\pi |x|})\biggr) \mathrm{d} w  \\
    =& \int_\R \frac{\mathscr{P}_y(w)}{p} \biggr(\log\left(\frac{\mathscr{P}_\alpha(x-w)}{(1-\alpha)\ee^{-\pi|x-w|}}\right)+ \log\left(\ee^{\pi (|x|-|x-w|)}\right)\biggr) \mathrm{d} w \\ 
    &+ \int_\R \frac{\mathscr{P}_{1-y}(w)}{q} \biggr( \log\left(\frac{\mathscr{P}_{1-\alpha}(x-w)}{\alpha \ee^{-\pi|x-w|}}\right) +\log\left(\ee^{\pi (|x|-|x-w|)}\right)\biggr) \mathrm{d} w\\
    \lesssim& \int_{\R} \left(\mathscr{P}_y(w)+ \mathscr{P}_{1-y}(w)\right) (1+|w|) \mathrm{d} w \lesssim 1,
\end{align*}
which completes the proof.
\end{proof}

We can now finish the proof of Theorem~\ref{thm:optimizers}.

\begin{proof}[Proof of Theorem~\ref{thm:optimizers}] Let $\alpha \in (0,1)$,$1\leq p,q\leq  \infty$. Then we define $h(z) = \ee^{\phi_{\alpha,p,q}(z)}$, where $\phi_{\alpha,p,q}$ is the function from Lemma~\ref{lem:phasefunction}. Then by~\eqref{eq:phasedecay} we have
\begin{align}
    \ee^{-\pi \left(\frac{1-y}{p} + \frac{y}{q}\right)|x|}\lesssim |h(x + \ii y)| \lesssim \ee^{-\pi \left(\frac{1-y}{p} + \frac{y}{q}\right)|x|},
\end{align}
uniformly in $0<y<1$; therefore,
\begin{align*}
    \norm{h_y}_{L^p(\R)} \lesssim (1-y)^{-\frac{1}{p}} \quad\mbox{and}\quad \norm{h_y}_{L^q(\R)} \lesssim y^{-\frac{1}{q}} \quad \mbox{for any $0<y<1$.} 
\end{align*}
In paticular, $h \in \mathbb{H}^{p,q}(S)$. Moreover, since
\begin{align*}
    \lim_{y \downarrow 0} \re \phi_{\alpha,p,q}(x+\ii y) = \frac{1}{p} \log\left(\frac{\mathscr{P}_\alpha}{1-\alpha}\right)(x) \quad\mbox{and}\quad \lim_{y \uparrow 1} \re \phi_{\alpha,p,q}(x+\ii y) = \frac{1}{q} \log\left(\frac{\mathscr{P}_{1-\alpha}}{\alpha}\right)(x)
\end{align*}
 for any $x\in \R$ by Lemma~\ref{lem:PoissonRepresentation}, we have 
\begin{align}
    |h_0(x)| = \biggr(\frac{\mathscr{P}_\alpha(x)}{1-\alpha}\biggr)^{\frac{1}{p}} \quad \mbox{and} \quad |h_1(x)| = \biggr(\frac{\mathscr{P}_{1-\alpha}(x)}{\alpha}\biggr)^{\frac{1}{q}} \quad \mbox{for almost every $x\in \R$.} \label{eq:valuesatboundary}
\end{align}
Next, let $\beta  \in \R$ be such that
\begin{align}
    \ee^{-\beta \alpha} h(\ii \alpha) = \frac{1}{1-\alpha} \label{eq:atitheta}
\end{align}
(which exists because $\mathrm{Im}\, \phi_{\alpha,p,q}(\ii \alpha) = 0$) and define $m: S \rightarrow \C$ via the formula
\begin{align*}
    m(z) \coloneqq  \frac{\mathscr{P}_\alpha(z)}{1-\alpha} \frac{\ee^{-\ii\beta z}}{h^\ast(z)} \quad z\in S^\ast. 
\end{align*}
 Thus, by \eqref{eq:valuesatboundary} and the definition of $m$, we see that $h \ee^{\ii \beta z}$ and $m$ satisfy the Euler-Lagrange equations~\eqref{eq:ELeq0} and \eqref{eq:ELeq1}. So by Lemma~\ref{lem:eulerlagrange}, it remains to show that $m \in \mathscr{P}_\alpha + \mathbb{H}^{p^\ast,q^\ast}(S)$. For this, we note that the poles of $m$ and $\mathscr{P}_\alpha$ at $-\ii \alpha$ cancel out (see~\eqref{eq:atitheta}) and therefore
\begin{align}
    |m(x+\ii y) - \mathscr{P}_\alpha(x+\ii y)| \lesssim 1 \quad \mbox{for any $x+\ii y$ in a neighborhood of $-\ii \alpha$.} \label{eq:boundedatitheta}
\end{align}
On the other hand, by~\eqref{eq:phasedecay} we find
\begin{align}
    |m(x+\ii y)| \lesssim \frac{|\mathscr{P}_\alpha(x+\ii y) \ee^{\pi |x|}|}{\ee^{\pi |x|} |h(x+\ii y)|} \lesssim \ee^{-\pi \left(\frac{1+y}{p^\ast} - \frac{y}{q^\ast}\right) |x|}, \quad \mbox{uniformly for $x+\ii y$ away from $-\ii \alpha$.} \label{eq:boundedawayitheta}
\end{align}
Combining~\eqref{eq:boundedatitheta} and \eqref{eq:boundedawayitheta}, we conclude that $m - \mathscr{P}_\alpha \in \mathbb{H}^{p^\ast,q^\ast}(S^\ast)$, which finishes the proof.
\end{proof}
\addtocontents{toc}{\protect\setcounter{tocdepth}{-1}}
\section*{Acknowledgements}
The author acknowledges funding by the \emph{Deutsche Forschungsgemeinschaft} (DFG, German Research Foundation) - Project number 442047500 through the Collaborative Research Center "Sparsity and Singular Structures" (SFB 1481).

\addtocontents{toc}{\protect\setcounter{tocdepth}{2}}
\appendix



\section{Measures with one-sided exponentially decaying Fourier transform}
\label{app:analyticmeasures}
In this section, we prove a regularity result for certain measures whose Fourier transform has an one-sided exponential decay in a tempered sense. Precisely, we prove the following result.

\begin{lemma}[Absolute continuity of measures with exponentially weighted Fourier transform]\label{lem:abscont} Let $\mu_0, \mu_1 \in \mathcal{M}(\R) + L^\infty(\R)$ and suppose that $\widehat{\mu}_1 = \exp \widehat{\mu_0} \in \mathcal{S}^\prime(\R)$. Then $\mu_0$ and $\mu_1$ are absolutely continuous with respect to the Lebesgue measure, i.e., $\mu_0,\mu_1 \in L^1(\R) + L^\infty(\R)$.
\end{lemma}

This result is a improved version of \cite[Lemma B.2]{CCR24}. Its proof relies on a classical result of F.\ and M.\ Riesz \cite{RR16} about \emph{analytic} measures, whose following version can be found in \cite[Lemma 13.4]{Mas09} (see also \cite[Theorem 3.8]{Gar06}).

\begin{lemma}[Lemma 13.4 \cite{Mas09}] \label{lem:Riesz} Let $\mu \in \mathcal{M}(\R)$. Then the following are equivalent:
\begin{enumerate}[label=(\roman*)]
\item $\mu$ is analytic;
\item for any $z\in \C_+ = \{z\in\CC: \mathrm{Im}(z) >0\}$, 
    \begin{align*}
	   \int_\R \frac{\mu (\mathrm{d}x)}{x-\overline{z}}  = 0; 
    \end{align*}
\item $\mu(\mathrm{d} x) = u(x) \mathrm{d}x$, where $u \in L^1(\R)$ with $\widehat{u}(k) = 0$ for all $k\leq 0$.
\end{enumerate}
\end{lemma}

\begin{proof}[Proof of Lemma~\ref{lem:abscont}] Let $\phi \in C^\infty(\R;[0,1])$ be a function satisfying $\phi(k) 
= 0$ for $k\leq -1$ and $\phi(k) =1$ for $k \geq 0$. We then define $\mu_+ \in \mathcal{S}^\prime$ via the formula
\begin{align*}
	\widehat{\mu}_+ = \widehat{\mu}_1 \exp_{-1} \phi = \widehat{\mu}_0 \phi.
 \end{align*}
 Note that, since $\psi = \exp_{-1} \phi \in \mathcal{S}(\R)$ (as the support of $\phi$ is contained in $(-1,\infty)$), we have
 \begin{align*}
     \mu_+ = \mu_0 \ast \psi \in L^\infty(\R) ,
 \end{align*}
 by Young's inequality. In particular
 \begin{align*}
     \mu_- \coloneqq \mu_0 - \mu_+ \in \mathcal{M}(\R) + L^\infty(\R) \quad \mbox{and} \quad \widehat{\mu}_-(k) = 0 \quad \mbox{for $k\leq 0$.}
 \end{align*}
Hence, it is enough to show that $\mu_-$ is absolutely continuous. For this, we first make the following observation: for any $f, g \in L^1(\R)$ with $f(x) = 0 = g(x)$ for any $x\leq 0$ we have $f \ast g \in L^1(\R)$ and $f\ast g(x) = 0$ for $x\leq 0$. In other words, the space of integrable functions with support on the non-negative real axis is a subalgebra of $L^1(\R)$. Consequently, for any $z_1, z_2 \in \C_+$, the function 
\begin{align*}
    f_{z_1,z_2}(x) = \frac{1}{(x-\overline{z}_1)^2}\frac{1}{x-\overline{z_2}} \quad \mbox{satisfies} \quad \widehat{f}_{z_1,z_2}(k) = 0 \quad \mbox{for $k\leq 0$.}
\end{align*}
Indeed, this follows from the previous observation and the fact that for any $z \in \C$, the Fourier transform of $f_z(x) = 1/(x-\overline{z})$ is given by $\widehat{f}_z(k) = \ee^{i \overline{z} k} \mathbb{1}_{k \geq 0}$, which is in $L^1(\R)$ and has support on the non-negative axis. Since the supports of $\widehat{\mu}_-$ and $\widehat{f}_{z_1,z_2}$ are disjoint, a standard approximation argument via mollifiers together with Plancherel/Parseval identity leads to the identity
\begin{align*}
    \int_{\R} f_{z_1,z_2}(x) \mu_-(\mathrm{d}x) = \int_{\R} \widehat{f}_{z_1,z_2}(k) \widehat{\mu}_-(k) \mathrm{d} k = 0 \quad \mbox{for any $z_1, z_2 \in \C_+$.}
\end{align*}
 But since $\mu_- \in \mathcal{M}(\R) + L^\infty(\R)$, the measure $\mu_-(d x)/(x-\overline{z_1})^2$ is a bounded Radon measure, and we conclude that $\mu_-$ is absolutely continuous by Lemma~\ref{lem:Riesz}. \end{proof}
\section{The Poisson kernel on the strip}
\label{app:Poisson}
In this section, we recall some results about the Poisson kernel on the strip,
\begin{align}
    \mathscr{P}_y(x) \coloneqq \frac{1}{2} \frac{\sin(\pi y)}{\cosh(\pi x) - \cos(\pi y)}, \quad (x,y) \in \R\times (0,1). \label{eq:Poissondef1}
\end{align}

The first result we recall about this kernel, and which justifies its name, is the following integral representation for harmonic functions on the strip. A proof can be found in \cite[Theorem 1]{Wid61}.

\begin{theorem}[Poisson integral on the strip]\label{lem:PoissonRepresentation} Let $f: \R \rightarrow \C$ and $g: \R \rightarrow \C$ be two measurable functions such that $\frac{f(t)}{\cosh(\pi t)}, \frac{g(t)}{\cosh(\pi t)} \in L^1(\R)$, then the function defined as
\begin{align*}
    u(x,y) \coloneqq \mathscr{P}_y \ast f + \mathscr{P}_{1-y} \ast g
\end{align*}
is harmonic in the strip $S$ and satisfies
\begin{align*}
    \lim_{(x,y) \ra (x_0,0)} u(x,y) = f(x_0) \quad\mbox{and}\quad \lim_{(x,y) \ra (x_0,1)} u(x,y) = g(x_0),
\end{align*}
at each continuity point $x_0 \in \R$.
\end{theorem}

The second result we use throughout the paper is an explicit formula for the Fourier transform of $\mathscr{P}_y$. Since we could not find a reference for this formula, we provide the simple proof below.

\begin{lemma}[Fourier transform of Poisson kernel]\label{lem:PoissonFourier} Let $\mathscr{P}_y$ be the Poisson kernel defined in \eqref{eq:Poissondef1}. Then we have
\begin{align}
    \widehat{\mathscr{P}}_y(k) = \int_\R \mathscr{P}_y(x) \ee^{-\ii k x} \mathrm{d} x = \frac{\sinh\left((1-y) k\right)}{\sinh(\pi k)}. \label{eq:PoissonFourier}
\end{align}
In particular, we have
\begin{align}
    \norm{\mathscr{P}_y}_{L^1} = \int_{\R} \mathscr{P}_y(x) \mathrm{d} x = 1-y. \label{eq:L1Poisson}
\end{align}
for any $0<y<1$.
\end{lemma}

\begin{proof}
    First, we observe that for any meromorphic function $f$ that is $2\ii$-periodic on the strip $S$, has no poles along the real axis, and decays fast enough as $|x| \ra \infty$, we have
    \begin{align}
        \widehat{f}(\omega) := \int_\R f(\tau) \ee^{-\imath\omega \tau} \mathrm{d} \tau = -\frac{\pi \imath}{\sinh(\omega)} \mathrm{Res}(\ee^{-\imath\omega(z-\imath)} f,S), \label{eq:FT1}
\end{align}
where $\mathrm{Res}(\ee^{-\imath\omega(z-\imath\pi)} f,S)$ denotes the residue of the function $f(z) \ee^{-\imath\omega(z-\imath\pi)}$ inside the strip $S$. Indeed, from the residue theorem applied to the contour $\partial S$ (which is possible by a standard limiting argument as long as $f$ decays fast enough as $|x| \ra \infty$),
\begin{align*}
    2\pi i \mathrm{Res}(f \ee^{-\ii \omega(z-\ii)}, S) =  \ee^{-\omega}\int_\R f(x) \ee^{-\ii\omega x} \mathrm{d} x - \ee^{\omega}\int_\R f_2(x) \ee^{-\ii \omega x} \mathrm{d} x = -2 \sinh(\omega) \widehat{f}(\omega),
\end{align*}
which proves~\eqref{eq:FT1}. 

Next, note that the Poisson kernel $z\in \C \mapsto \mathscr{P}_y(z)$ is a $2\ii$-periodic merormophic function in $\C$ with single poles along the set $\{ \pm \ii y + 2\ii k: k \in \Z\}$. Moreover, it decays (exponentially) fast as $|x| \ra\infty$. Formula~\eqref{eq:PoissonFourier} now follows from~\eqref{eq:FT1} and straightforward calculations.

Equation~\eqref{eq:L1Poisson} follows by noticing that $\mathscr{P}_y(x) > 0$ for $x \in \R$ and taking the limit $\lim_{\omega \ra 0} \widehat{\mathscr{P}}_y(\omega)$.
\end{proof}
\section{Analytic formula for optimal values}
\label{app:analyticformula}
In this section we prove the formula from Remark~\ref{rem:analyticformula}. To this end, we shall use the following properties of the dilogarithm function
\begin{align}
    \mathrm{Li}_2(z) \coloneqq - \int_0^z \frac{\log(1-u)}{u} \mathrm{d} u, \quad z\in \C\setminus [1,\infty). \label{eq:dilog}
\end{align}
(For a proof, see \cite[Equations (1.7), (1.11), and (1.15)]{Lew81}.)
\begin{lemma}[Dilogarithm properties] \label{lem:dilogproperties} Let $\mathrm{Li}_2:\C\setminus [1,\infty)$ be the function defined in \eqref{eq:dilog}. Then $\mathrm{Li}_2$ is holomorphic and satisfies the identities
\begin{align}
    &\mathrm{Li}_2(z) + \mathrm{Li}_2(1/z) = -\mathrm{Li}_2(1) - \frac{\log(-z)^2}{2} \quad &\mbox{for any $z\in \C \setminus [0, \infty)$,} \label{eq:dilogid0} \\
    &\mathrm{Li}_2(1-z) + \mathrm{Li}_2(z) = \mathrm{Li}_2(1)-\log(z) \log(1-z), \quad &z \not \in (-\infty,0] \cup [1, \infty), \label{eq:dilogid1} \\
    &\frac{1}{2} \mathrm{Li}_2(z^2) = \mathrm{Li}_2(-z) + \mathrm{Li}_2(z), \quad &z\not \in (-\infty,-1) \cup (1,\infty),\label{eq:dilogid2}
\end{align}
\end{lemma}

\begin{proof}[Proof of formula~\eqref{eq:analyticformula}]
First, we recall that the Poisson kernel on the strip $S$ is given by
\begin{align*}
    \mathscr{P}_{\alpha}(x) = \frac{1}{2} \frac{\sin(\pi\alpha)}{\cosh(\pi x) - \cos(\pi \alpha)}
\end{align*}
Thus using that $\int_\R \mathscr{P}_{\alpha}(x) \mathrm{d} x = 1-\alpha$ (see Lemma~\ref{lem:PoissonFourier}) in \eqref{eq:H1infvalue} we find
\begin{align}
    \log H_{1,\infty}(\alpha) &= \int_{\R} \mathscr{P}_{\alpha}(\tau) \log\left(\frac{\mathscr{P}_{\alpha}(x)}{1-\alpha}\right)\mathrm{d} x \nonumber \\
    &=\int_\R  \mathscr{P}_{\alpha}(x) \log\left(\frac{\sin(\pi \alpha)}{2 (1-\alpha)}\right)\mathrm{d} x - \frac{\sin(\pi \alpha)}{2} \int_\R \frac{\log\left( \cosh(\pi x) - \cos(\pi \alpha)\right)}{\cosh(\pi x) - \cos(\pi \alpha)} \mathrm{d} x  \nonumber \\
    &= (1-\alpha)\log\left(\frac{\sin(\pi \alpha)}{2(1-\alpha)}\right) - \underbrace{\sin(\pi \alpha)\int_{-\infty}^0 \frac{\log\left(\cosh(\pi x) - \cos(\pi \alpha)\right)}{\cosh(\pi x) - \cos(\pi \alpha)} \mathrm{d} x}_{\coloneqq I(\alpha)}. \label{eq:splitingintegrals0}
\end{align}
To evaluate the integral $I(\alpha)$ we first use the change of variables $t = \e^{\pi x}$, $\mathrm{d} x = \mathrm{d} t/(\pi t)$ to find
\begin{align}
    I(\alpha) &= \frac{2 \sin(\pi \alpha)}{\pi} \int_0^1 \frac{\log\left(t^2 - 2t \cos(\pi \alpha) + 1\right)- \log(2t)}{t^2 - 2t \cos(\pi \alpha) + 1} \mathrm{d} t \nonumber \\
    &= \underbrace{\frac{2 \sin(\pi \alpha)}{\pi}\int_0^1 \frac{\log(t^2-2t \cos(\pi \alpha) +1)}{t^2-2t \cos(\pi \alpha) + 1} \mathrm{d} t}_{\coloneqq I_1(\alpha)} - \underbrace{\frac{2 \sin(\pi \alpha)}{\pi}\int_0^1 \frac{\log(2t)}{t^2 - 2t \cos(\pi \alpha) + 1} \mathrm{d} t}_{\coloneqq I_2(\alpha)}.  \label{eq:splittingintegrals1}
\end{align}
We now note that $t^2-2t \cos(\pi \alpha) +1 = (t - \ee^{\ii \pi \alpha})(t-\ee^{-\ii \pi \alpha})$, and therefore,
\begin{align*}
    \frac{2\sin(\pi \alpha)}{t^2-2t \cos(\pi \alpha) +1} = \frac{1}{\ii } \biggr( \frac{1}{t-\ee^{\ii \pi \alpha}}- \frac{1}{t-\ee^{-\ii \pi \alpha}}\biggr).
\end{align*}
So from the change of variables $t = \ee^{\ii \pi \alpha} - \ee^{\pi \alpha} \tau$ and \eqref{eq:dilogid2}, the term $I_2(\alpha)$ can be written as
\begin{align}
    I_2(\alpha) &= \frac{1}{\ii \pi} \biggr(\int_0^1 \frac{\log(t)+\log(2)}{t - \ee^{\ii \pi \alpha}} \mathrm{d}t - \int_0^1 \frac{\log(t)+\log(2)}{t-\ee^{-\ii \pi \alpha}} \mathrm{d} t\biggr)\nonumber\\
    &= \frac{2}{\pi} \mathrm{Im} \, \biggr( \int_{1}^{1-\ee^{-\ii \pi \alpha}} \frac{\log\left(\ee^{\ii \pi \alpha}(1-\tau)\right)}{\tau} \mathrm{d} \tau + \log(2)\left(\log(1-\ee^{\ii \pi \alpha})-\log(-\ee^{\ii \pi \alpha}\right)\biggr) \nonumber \\
    &= \frac{2}{\pi} \mathrm{Im}\biggr(\mathrm{Li}_2(1) - \mathrm{Li}_2\left(1-\ee^{-\ii \pi \alpha}\right) - \log\left(\ee^{-\ii \pi \alpha}\right)\log\left(1-\ee^{-\ii \pi \alpha}\right) \biggr)+ (1-\alpha)\log(2)  \nonumber \\ 
    &\stackrel{\eqref{eq:dilogid2}}{=} \frac{2}{\pi} \mathrm{Im} \,\mathrm{Li}_2(\ee^{-\ii \pi \alpha}) + (1-\alpha) \log(2) . \label{eq:I2result}
\end{align}
For $I_1(\alpha)$ we use the same splitting to obtain
\begin{align*}
    I_1(\alpha) &= \frac{2 \sin(\pi \alpha)}{\pi} \int_0^1 \frac{\log(t-\ee^{\ii \pi \alpha}) + \log(t-\ee^{-\ii \pi \alpha})}{t^2 - 2t \cos(\pi \alpha) +1} \mathrm{d} t \\
    &= \underbrace{\frac{1}{\ii \pi}\int_0^1 \frac{\log(t-\ee^{\ii \pi \alpha})}{t-\ee^{\ii \pi \alpha}} - \frac{\log(t-\ee^{-\ii \pi \alpha})}{t-\ee^{-\ii \pi \alpha}} \mathrm{d} t}_{\coloneqq I_{1,1}(\alpha)} + \underbrace{\frac{1}{\ii \pi}\int_0^1 \frac{\log(t-\ee^{-\ii \pi \alpha})}{t-\ee^{\ii \pi \alpha}} - \frac{\log(t-\ee^{\ii \pi \alpha})}{t-\ee^{-\ii \pi \alpha}} \mathrm{d} t}_{\coloneqq I_{1,2}(\alpha)}. 
\end{align*}
The first part can be directly integrated using the substitution $u= t-\ee^{\pm \ii \pi \alpha}$ to obtain
\begin{align}
    I_{1,1}(\alpha) &= \frac{\log(t-\ee^{\ii \pi \alpha})^2-\log(t-\ee^{-\ii\pi \alpha})^2}{2 \pi \ii}\biggr\rvert_{t=0}^{t=1}= - \frac{1-\alpha}{2} \log\left(2\left(1-\cos(\pi\alpha)\right)\right) ,\label{eq:I21result}
\end{align}
where we used the trigonometric identity $\arctan(\cot(\pi\alpha/2)) = \pi/2-\pi\alpha/2$ (which holds for any $0\leq\alpha\leq 1$) in the last step. To evaluate $I_{1,2}(\alpha)$ we note that
\begin{align*}
    I_{1,2}(\alpha) = \frac{1}{\pi \ii} \int_0^1 \frac{\log(t-\ee^{-\ii \pi \alpha})}{t-\ee^{\ii \pi \alpha}} - \overline{\frac{\log(t-\ee^{-\ii \pi \alpha})}{t-\ee^{\ii \pi \alpha}}} \mathrm{d} t = \frac{2}{\pi} \mathrm{Im}\, \int_0^1 \frac{\log(t-\ee^{-\ii \pi \alpha})}{t-\ee^{\ii \pi \alpha}} \mathrm{d} t
\end{align*}
and therefore it is enough to evaluate the integral $\int_0^1 \frac{\log(t-\ee^{-\ii \pi \alpha})}{t-\ee^{\ii \pi \alpha}} \mathrm{d} t$. For this, we first make the changes of variable $\tau = -t+\ee^{\ii \pi \alpha}$ and $s = \tau/(2\ii \sin(\pi \alpha))$ to obtain
\begin{align}
    \int_0^1 \frac{\log(t-\ee^{-\ii \pi \alpha})}{t-\ee^{\ii \pi \alpha}} \mathrm{d}t =& \int_{\ee^{\ii \pi \alpha}}^{\ee^{\ii \pi \alpha}-1} \frac{\log(2 \ii \sin(\pi \alpha) - \tau)}{\tau} \mathrm{d} \tau \nonumber \\
    =& \int_{\frac{\ee^{\ii \pi \alpha}}{2\ii \sin(\pi \alpha)}}^{\frac{\ee^{\ii \pi \alpha}-1}{2\ii \sin(\pi \alpha)}} \frac{\log\left(2 \ii \sin(\pi \alpha)\right) + \log(1-s)}{s} \mathrm{d} t.\nonumber \\
    =& \log\left(2\ii \sin(\pi \alpha)\right)\left(\log\left(\frac{\ee^{\ii \pi \alpha}-1}{2\ii \sin(\pi \alpha)}\right)-\log\left(\frac{\ee^{\ii \pi \alpha}}{2\ii \sin(\pi \alpha)}\right)\right) \nonumber \\
    &+ \mathrm{Li}_2\left(\frac{\ee^{\ii \pi \alpha}}{2\ii \sin(\pi \alpha)}\right) - \mathrm{Li}_2\left(\frac{\ee^{\ii \pi \alpha}-1}{2\ii \sin(\pi \alpha)}\right). \label{eq:I12middle}
\end{align}
Hence, using that $\log\left(2\ii \sin(\pi \alpha)\right) = \log\left(2\sin(\pi \alpha)\right)+ \ii \frac{\pi}{2}$ and $\log\left(\frac{\ee^{\ii \pi \alpha}-1}{2\ii \sin(\pi \alpha)}\right)-\log\left(\frac{\ee^{\ii \pi \alpha}}{2\ii \sin(\pi \alpha)}\right) 
    = \log\left(2\sin(\pi\alpha/2)\right) + \ii \frac{\pi}{2} - \ii \frac{\pi \alpha}{2}$ in \eqref{eq:I12middle}, we find
\begin{align}
I_{1,2}(\alpha) = (1-\alpha) 
 \log\left(2\sin(\pi \alpha)\right) + \frac{1}{2} \log&\left(2 \sin\left(\frac{\pi \alpha}{2}\right)\right)\nonumber \\ &+ \frac{2}{\pi} \mathrm{Im}\,\left(\mathrm{Li}_2\left(\frac{\ee^{\ii \pi \alpha}}{2\ii \sin(\pi \alpha)}\right) - \mathrm{Li}_2\left(\frac{\ee^{\ii\pi \alpha}-1}{2\ii \sin(\pi \alpha)}\right)\right). \label{eq:I12result}
\end{align}
Putting~\eqref{eq:I2result}, \eqref{eq:I21result}, and \eqref{eq:I12result} back in \eqref{eq:splittingintegrals1} and substituting the result in \eqref{eq:splitingintegrals0} we find
\begin{align}
    \log H_{1,\infty}(\alpha) =& (1-\alpha)\log\left(\frac{\sin(\pi \alpha)}{2(1-\alpha)}\right) - I_{1,1}(\alpha)-I_{1,2}(\alpha)+I_2(\alpha) \nonumber \\
    =& (1-\alpha)\log\left(\frac{\sin(\pi \alpha)}{2(1-\alpha)}\right) + \frac{1-\alpha}{2} \log\left(2\left(1-\cos(\pi \alpha)\right)\right) - (1-\alpha) \log\left(2\sin(\pi \alpha)\right) \nonumber \\
    &- \frac{1}{2} \log\left(2\left(1-\cos(\pi \alpha)\right)\right) + (1-\alpha) \log(2) \nonumber \\
    &-\frac{2}{\pi} \mathrm{Im}\,\biggr(\mathrm{Li}_2\left(\frac{\ee^{\ii \pi \alpha}}{2\ii \sin(\pi \alpha)}\right) - \mathrm{Li}_2\left(\frac{\ee^{\ii\pi \alpha}-1}{2\ii \sin(\pi \alpha)}\right) - \mathrm{Li}_2(\ee^{-\ii \pi \alpha})\biggr)\nonumber \\
    =& - \log\left(2 (1-\alpha)^{1-\alpha} \sin\left(\frac{\pi \alpha}{2}\right)^\alpha\right)  -\underbrace{\frac{2}{\pi} \mathrm{Im}\biggr(\mathrm{Li}_2\left(\frac{\ee^{\ii \pi \alpha}}{2\ii \sin(\pi \alpha)}\right) - \mathrm{Li}_2\left(\frac{\ee^{\ii\pi \alpha}-1}{2\ii \sin(\pi \alpha)}\right) - \mathrm{Li}_2(\ee^{-\ii \pi \alpha}) \biggr)}_{\coloneqq T(\alpha)}.\label{eq:fast}
\end{align}
Finally, by setting $z= \ee^{-\ii \pi \alpha}$ we note that by Lemma~\ref{lem:dilogproperties},
\begin{align*}
    T(\alpha) =& \frac{2}{\pi}\mathrm{Im}\,\biggr(\mathrm{Li}_2\left(\frac{1}{1-z^2}\right) - \mathrm{Li}_2\left(\frac{1}{1+z}\right) - \mathrm{Li}_2(z)\biggr), \\
    \stackrel{\eqref{eq:dilogid0}}{=}& \frac{2}{\pi}\mathrm{Im}\,\biggr(-\mathrm{Li}_2(1-z^2) +\mathrm{Li}_2(1+z)  - \mathrm{Li}_2(z)-\frac{\log(z^2-1)^2}{2}+\frac{\log(-1-z)^2}{2}\biggr)\\
    \stackrel{\eqref{eq:dilogid1}}{=}& \frac{2}{\pi}\mathrm{Im}\,\biggr(\mathrm{Li}_2(z^2) -\mathrm{Li}_2(-z) - \mathrm{Li}_2(z) +\log(z^2)\log(1-z^2)- \log(-z) \log(1+z) -\frac{\log(z^2-1)^2}{2} \\  &+ \frac{\log(-1-z)^2}{2} \biggr) \\
    \stackrel{\eqref{eq:dilogid2}}{=}& \frac{2}{\pi} \mathrm{Im}\, \biggr(\frac{1}{2} \mathrm{Li}_2(z^2)+\log(z^2)\log(1-z^2)- \log(-z) \log(1+z) -\frac{\log(z^2-1)^2}{2} + \frac{\log(-1-z)^2}{2}\biggr)  \nonumber \\
    &= \frac{1}{2} \mathrm{Im} \, \mathrm{Li}_2(\ee^{-\ii 2\pi \alpha}) - \log\left(\frac{\sin\left(\frac{\pi \alpha}{2}\right)^\alpha}{2^{1-2\alpha} \sin(\pi \alpha)^{1-\alpha}}\right), \quad 
\end{align*}
and therefore, plugging the above expression back in \eqref{eq:fast}, we obtain
\begin{align*}
    \log H_{1,\infty}(\alpha) =  -(1-\alpha) \log\left(4 (1-\alpha) \sin(\pi \alpha)\right) + \frac{1}{\pi} \mathrm{Im} \,\mathrm{Li}_2(\ee^{\ii 2\pi \alpha}), 
\end{align*}
which completes the proof. 
\end{proof}
\section{Weigthed three-lines lemma}
\label{app:weigthed}
We now prove Theorem~\ref{thm:weightedversion}.

\begin{proof}[Proof of Theorem~\ref{thm:weightedversion}]
From the assumptions on $w_0$ and $w_1$, the functions $\log w_0:\R \rightarrow \R$ and $\log: \R \rightarrow\R$ are bounded and continuous. Therefore, by Lemma~\ref{lem:PoissonRepresentation} the function
\begin{align*}
    u(x,y) = (P_y \ast \log w_0)(x) + (P_{1-y}\ast \log w_1)(x)
\end{align*}
is the unique bounded harmonic function in $\mathcal{S}$ satisfying 
\begin{align}
    \lim_{y \downarrow 0} u(x,y) = \log w_0(x)\quad \mbox{and}\quad \lim_{y \uparrow 1} u(x,y) = \log w_1(x), \quad \mbox{for any $x \in \R$.} \label{eq:boundaryconditions}
\end{align} 
Hence, there exists a unique (up to a constant) harmonic conjugate on $\mathcal{S}$, which we denote by $v$, and the function
\begin{align*}
    z = x+\ii y \mapsto f(x+\ii y) \coloneqq \exp\left(u(x,y) + \ii v(x, y)\right)
\end{align*}
is holomorphic and uniformly bounded in $\mathcal{S}$. Moreover, as $f_y \rightarrow f_0$ and $f_y \rightarrow f_1$ a.e. when $y \downarrow 0$, respectively $y \uparrow 1$ for some $f_0,f_1 \in \mL^\infty(\R)$, by~\eqref{eq:boundaryconditions}, we see that $|f_y| \rightarrow w_0$ almost everywhere. Consequently, for any $h\in \mathbb{H}^{p,q}(\mathcal{S})$, we have $g \coloneqq h f \in \mathbb{H}^{p,q}(\mathcal{S})$ and 
\begin{align*}
    \norm{h_0 w_0}_{\mL^p} = \norm{g_0}_{\mL^p} \quad \mbox{and}\quad \norm{h_1 w_1}_{\mL^q} = \norm{g_1}_{\mL^q}.
\end{align*}
The result now follows by applying Theorem~\ref{thm:optimum} to $g$. The sharpness follows by choosing $h$ such that $g$ is the optimizer of~\eqref{eq:pq3lineproblem}.
\end{proof}

\section*{Data availability}
No datasets were generated or analysed during the current study.

\bigskip

\begin{thebibliography}{EFGHW21}

\bibitem[BK07]{BK07}
	\textsc{A.~Bakan} and \textsc{S.~Kaijser},
	\newblock Hardy spaces for the strip.
	\newblock \doi{10.1016/j.jmaa.2006.10.088}{\emph{Journal of Mathematical Analysis and Applications}} \textbf{333} (2007), no.1, 347--364.
	\newblock \mr{2323495}.
	\newblock \zbl{1160.30021}.
	\hfill

 \bibitem[CCR24]{CCR24}
	\textsc{T.C.~Corso} and \textsc{T. Ried},
	\newblock {\emph{On a variational problem related to the Cwikel-Lieb-Rozenblum and Lieb-Thirring inequalities}}.
	\newblock \href{https://arxiv.org/abs/2403.04347}{arxiv 2403.04347}
	\hfill

\bibitem[CKP07]{CKP07}
	\textsc{R.D.~Carmichael, A.~Kami\'nski}, and \textsc{S.~Pilipovi\'c},
	\newblock \doi{10.1142/9789812708786}{\emph{Boundary values and convolution in ultradistribution spaces}}.
	\newblock  Series on Analysis, Applications and Computation Vol.\ 1, World Scientific, Hackensack, NJ, 2007.
	\newblock \mr{2347838}.
	\newblock \zbl{1155.46002}.
	\hfill


\bibitem[FHJN21]{FHJN21}
	\textsc{R.L.~Frank, D.~Hundertmark, M.~Jex}, and \textsc{P.T.~Nam},
	\newblock The Lieb--Thirring inequality revisited.
	\newblock \doi{10.4171/jems/1062}{\emph{Journal of the European Mathematical Society (JEMS)}} \textbf{23} (2021), no. 8, 2583--2600. 
	\newblock \mr{4269422}.
	\newblock \zbl{1467.35014}.
	\hfill

\bibitem[FLW23]{FLW23}
	\textsc{R.L.~Frank, A.~Laptev}, and \textsc{T.~Weidl},
	\newblock \doi{10.1017/9781009218436}{\emph{Schr\"odinger operators: eigenvalues and Lieb--Thirring inequalities}}.
	\newblock Cambridge Studies in Advanced Mathematics Vol.\ 200, Cambridge University Press, Cambridge, 2023.
	\newblock \mr{4496335}.
	\newblock \zbl{07595814}.
	\hfill

\bibitem[Fra14]{Fra14}
	\textsc{R.L.~Frank}
	\newblock Cwikel's theorem and the CLR inequality.
	\newblock \doi{10.4171/JST/59}{\emph{Journal of Spectral Theory}} \textbf{4} (2014), no.1, 1--21.
	\newblock \mr{3181383}.
	\newblock \zbl{1295.35347}.
	\hfill
 
\bibitem[Fra21]{Fra21}
	\textsc{R.L.~Frank},
	\newblock \doi{10.1090/pspum/104/01877}{The Lieb--Thirring inequalities: recent results and open problems}.
	\newblock In \emph{Nine mathematical challenges—an elucidation} (eds.\ A.~Kechris, N.~Makarov, D.~Ramakrishnan, X.~Zhu),
	\newblock  Proceedings of Symposia in Pure Mathematics \textbf{104}, American Mathematical Society, Providence, RI, 2021, pp. 210--235.
	\newblock \mr{4337417}.
	\newblock \zbl{1518.35002}.
	\hfill
 
\bibitem[Gar06]{Gar06}
	\textsc{J.B.~Garnett},
	\newblock \doi{10.1007/0-387-49763-3}{\emph{Bounded analytic functions}}. Revised 1st ed.
	\newblock Graduate Texts in Mathematics Vol.\ 236. Springer, New York, NY, 2006.
	\newblock \mr{2261424}.
	\newblock \zbl{1106.30001}.
	\hfill

 \bibitem[Gra14]{Gra14}
	\textsc{L.~Grafakos},
	\newblock \doi{10.1007/0-387-49763-3}{\emph{Classical Fourier Analysis}}. third ed.
	\newblock Graduate Texts in Mathematics Vol.\ 249. Springer, New York, NY, 2014.
	\newblock \mr{3243734}.
	\newblock \zbl{1304.42001}.
	\hfill



\bibitem[HKRV23]{HKRV23}
	\textsc{D.~Hundertmark, P.~Kunstmann, T.~Ried}, and \textsc{S.~Vugalter},
	\newblock Cwikel's bound reloaded.
	\newblock \doi{10.1007/s00222-022-01144-7}{\emph{Inventiones Mathematicae}} \textbf{231} (2023), 111--167.
	\newblock \mr{4526822}.
	\newblock \zbl{1510.35193}.
	\hfill
	


\bibitem[Lan99]{Lan99}
	\textsc{Serge.~Lang},
	\newblock \doi{10.1007/978-1-4757-3083-8}{\emph{Complex analysis}}. fourth edition. 
	\newblock  Graduate Texts in Mathematics Vol.\ 103, Springer-Verlag, New York, 1999.
	\newblock \mr{1659317}.
	\newblock \zbl{0933.30001}.
	\hfill



\bibitem[LW00]{LW00}
	\textsc{A.~Laptev} and \textsc{T.~Weidl},
	\newblock Sharp Lieb-Thirring inequalities in high dimensions. 
	\newblock \doi{10.1007/BF02392782}{\emph{Acta Mathematica}} \textbf{184} (2000), 87--111. 
	\newblock \mr{1756570}
	\newblock \zbl{1142.35531}
	\hfill


\bibitem[Lew81]{Lew81}
	\textsc{L.~Lewin},
	\newblock Polylogarithms and associated functions
	\newblock North-Holland Publishing, New York-Amsterdam, 1981
	\newblock \mr{618278}
	\newblock \zbl{0465.33001}
	\hfill

\bibitem[Lie76]{Lie76} 
	\textsc{E.H.\ Lieb}, 
	\newblock Bounds on the eigenvalues of the Laplace and Schr\"odinger operators. 
	\newblock \doi{10.1090/S0002-9904-1976-14149-3}{\emph{Bulletin of the American Mathematical Society}} \textbf{82} (1976), 751--753. 
	\newblock \mr{0407909}.
	\newblock \zbl{0329.35018}.
	\hfill 

	
\bibitem[Mas09]{Mas09}
	\textsc{J.~Mashreghi},
	\newblock \doi{10.1017/CBO9780511814525}{\emph{Representation theorems in Hardy spaces}}. 
	\newblock London Mathematical Society Student Texts Vol.\ 74, Cambridge University Press, Cambridge, 2009.
	\newblock \mr{2500010}.
	\newblock \zbl{1169.22001}.
	\hfill

\bibitem[Nam21]{Nam21}
	\textsc{P.T.~Nam},
	\newblock Direct methods to Lieb--Thirring kinetic inequalities.
	\newblock In \emph{Density functionals for many-particle systems. Mathematical theory and physical applications of effective equations} (eds.\ B.-G.~Englert, H.~Siedentop, M.-I.~Trappe).
	\newblock Lecture Notes Series, Institute for Mathematical Sciences, National University of Singapore Vol.\ 41, World Scientific, Hackensack, NJ, 2023, 81–115.
	\newblock \mr{4635285}.
	\newblock \zbl{1515.81017}.
	\hfill
	
\bibitem[Roc72]{Roc72}
	\textsc{R.T.~Rockafellar},
	\newblock \emph{Convex Analysis}.
	\newblock Princeton Mathematical Series Vol.\ 28, Princeton University Press, Princeton, NJ, 1970.
	\newblock \mr{0274683}.
	\newblock \zbl{0193.18401}.
	\hfill
	
\bibitem[RR16]{RR16}
	\textsc{F.~Riesz} and \textsc{M.~Riesz},
	\newblock Über die Randwerte einer analytischen Funktion.
	\newblock \emph{Quatrième congrès des math.\ scand.} (1916), 27--44.
	\newblock \zbl{47.0295.03}.
	\hfill
 
\bibitem[RS75]{RS75}
	\textsc{M.~Reed} and \textsc{B.~Simon},
	\newblock \emph{Methods of modern mathematical physics II: Fourier Analysis, self-adjointness}.
	\newblock Academic Press, Princeton, New York, 1975.
	\newblock \zbl{0308.47002}.
	\hfill
 
\bibitem[Rud87]{Rud87}
	\textsc{W.~Rudin},
	\newblock \emph{Real and complex analysis}.
	\newblock McGraw-Hill, New York, third edition, 1987
	\newblock \mr{924157}.
	\newblock \zbl{0925.00005}.
	\hfill
 

 \bibitem[Rum10]{Rum10}
	\textsc{M.~Rumin},
        \newblock \emph{Spectral density and Sobolev inequalities for pure and mixed states}.
	\newblock \doi{10.1007/s00039-010-0075-6}{\emph{Geometric and Functional Analysis}}, \textbf{20} (2010), no.3, 817--844.
	\newblock \mr{2720233}.
	\newblock \zbl{1218.58021}.
	\hfill
	
\bibitem[Sch22]{Sch22}
	\textsc{L.~Schimmer},
	\newblock \doi{10.4171/90-2/39}{The state of the Lieb--Thirring conjecture}.
	\newblock In \emph{The physics and mathematics of Elliott Lieb—the 90th anniversary. Vol. II} (eds.\ R.L.~Frank, A.~Laptev, M.~Lewin, R.~Seiringer), EMS Press, Berlin, 2022, 253--275.
	\newblock \mr{4531363}.
	\newblock \zbl{1500.81037}.
	\hfill
	
\bibitem[Til61]{Til61}
	\textsc{H.G.~Tillmann},
	\newblock Darstellung der Schwartzschen Distributionen durch analytische Funktionen.
	\newblock \doi{10.1007/BF01180167}{\emph{Mathematische Zeitschrift}} \textbf{77} (1961), 106--124.
	\newblock \mr{0139002}.
	\newblock \zbl{0099.09703}.
	\hfill

 \bibitem[Wid61]{Wid61}
	\textsc{D.V.~Widder},
	\newblock Functions harmonic in a strip.
	\newblock \doi{10.2307/2034126}{\emph{Proceedings of the American Mathematical Society}} \textbf{12} (1961), no. 1, 67--72.
	\newblock \mr{132838}.
	\newblock \zbl{0096.07703}.
	\hfill
	
		
\end{thebibliography}
\end{document}